\documentclass[12pt]{amsart}
\usepackage{mathrsfs}
\usepackage{amssymb}
\usepackage{amsbsy}
\usepackage{amsmath}
\usepackage{mathtools,booktabs}
\usepackage{graphicx}
\usepackage{caption}
\usepackage{subfig}
\usepackage[colorlinks,linkcolor=blue,citecolor=blue,pagebackref,hyperfootnotes=false]{hyperref}

\makeatletter
\renewcommand{\subsection}{\@startsection{subsection}{2}{\z@}%
  {.5\linespacing\@plus.7\linespacing}{.25\linespacing}%
  {\normalfont\normalsize\bfseries\raggedright}}
\renewcommand{\subsubsection}{\@startsection{subsubsection}{3}{\z@}%
  {.4\linespacing\@plus.3\linespacing}{.15\linespacing}%
  {\normalfont\normalsize\itshape\raggedright}}
\makeatother

\newtheorem{theorem}{Theorem}[section]
\newtheorem{lemma}[theorem]{Lemma}

\theoremstyle{definition}

\theoremstyle{remark}
\newtheorem{remark}[theorem]{Remark}

\allowdisplaybreaks
\numberwithin{equation}{section}

\newcommand{\R}{\mathbb R}
\newcommand{\Sph}{\mathbb S^{n-1}}
\newcommand{\Kn}{\mathcal K^n}
\newcommand{\Ko}{\mathcal K_o^n}
\newcommand{\Ke}{\mathcal K_e^n}

\newcommand{\ind}{\mathbf 1}

\newcommand{\comb}[3]{{#1}^{1-#3}{#2}^{#3}}
\DeclareMathOperator{\intt}{int}
\DeclareMathOperator{\supp}{supp}

\begin{document}
\begin{center}
{\Large\bfseries A unified approach to logarithmic Brunn--Minkowski inequalities and logarithmic Minkowski problems for geometric functionals}
\end{center}

\vskip 20pt

\begin{center}
{{\bfseries Yude~~Liu\textsuperscript{1}\quad\quad Xinbao~~Lu\textsuperscript{2}\quad\quad Qiang~~Sun\textsuperscript{3}}\\[6pt]
\small{\textsuperscript{1}School of Mathematical Sciences, Fudan University, Shanghai, 200433, China\\
\textsuperscript{2}School of Mathematical Sciences, Key Laboratory of Intelligent Computing and Applications (Tongji University), Ministry of Education, Tongji University, Shanghai, 200092, China\\
\textsuperscript{3}School of Mathematics and Statistics, Xi'an Jiaotong University, Xi'an, 710049, China}}
\end{center}

\vskip 10pt

\footnotetext{E-mail addresses: 1.\ liuyude@fudan.edu.cn; 2.\ xinbaolu@tongji.edu.cn; 3.\ sunqiang@fudan.edu.cn.}
\footnotetext{Research of the first author is supported by China Postdoctoral Science Foundation No.~2025M773080. Research of the second author is supported by NSFC No.~12501078. Research of the third author is supported by China Postdoctoral Science Foundation No. 2025M783084 and China National Postdoctoral Program for Innovative Talents No.~BX20250071.}

\begin{center}
\begin{minipage}{15cm}
{\bfseries Abstract} We establish two general principles for geometric functionals, one for proving logarithmic Brunn–Minkowski inequalities and the other for solving the associated logarithmic Minkowski problems. For origin-symmetric convex bodies, we prove logarithmic Brunn–Minkowski inequalities for the first Dirichlet eigenvalue, spectral zeta functions, torsional rigidity, and chord integrals, along with their equality conditions. For each functional, we give the necessary and sufficient condition for the existence of an origin-symmetric solution to the corresponding logarithmic Minkowski problem for even measures and prove uniqueness.

\end{minipage}
\end{center}

\vskip 25pt
\section{\bf Introduction}
\vskip 5pt

The setting of this article is the $n$-dimensional Euclidean space $\R^n$, $n\geqslant2$. A \emph{convex body} in $\R^n$ is a compact convex set with nonempty interior. Denote by $\Kn$, $\Ko$, and $\Ke$ the classes of convex bodies, convex bodies containing the origin in their interiors, and origin-symmetric convex bodies, respectively.

The Brunn--Minkowski theory of convex bodies has two central topics: the isoperimetric problem and the Minkowski problem. The classical isoperimetric inequality has a far-reaching generalization in the Brunn--Minkowski inequality, a cornerstone of convex geometry linking two basic notions: vector addition and volume. It asserts that for $K,L\in\Kn$ and $0<\theta<1$, 
\[
 |(1-\theta)K+\theta L|^{\frac{1}{n}} \geqslant (1-\theta)|K|^{\frac{1}{n}}+\theta|L|^{\frac{1}{n}},
\]
with equality if and only if $K$ and $L$ are homothetic. Here $|K|$ denotes the $n$-dimensional volume (i.e., Lebesgue measure) of $K$. By homogeneity, the Brunn--Minkowski inequality is equivalent to $|(1-\theta)K+\theta L|\geqslant |K|^{1-\theta}|L|^\theta$.

In 2012, B\"or\"oczky--Lutwak--Yang--Zhang \cite{BLYZ12} posed the \emph{logarithmic Brunn--Minkowski inequality conjecture}, which states that for $K,L\in\Ke$ and $0<\theta<1$,
\[
 |\comb K L\theta|\geqslant |K|^{1-\theta}|L|^\theta,
\]
where $\comb K L\theta$ is the logarithmic Minkowski combination of $K$ and $L$ defined by
\[
 \comb K L\theta=\{x\in\R^n:x\cdot u\leqslant h_K(u)^{1-\theta}h_L(u)^\theta\text{ for }u\in\Sph\}.
\]
Here $h_K(u)=\max\{x\cdot u:x\in K\}$ is the support function of $K$. Applying the arithmetic--geometric mean inequality to $h_K$ and $h_L$ yields $\comb K L\theta\subseteq(1-\theta)K+\theta L$. Thus the logarithmic inequality is \emph{stronger} than the classical one for origin-symmetric bodies. 

B\"or\"oczky--Lutwak--Yang--Zhang \cite{BLYZ12} themselves resolved the planar case. Saroglou \cite{SarUnconditional} proved the unconditional case. He later showed that the logarithmic Brunn--Minkowski inequality for volume implies the same inequality for every even log-concave measure \cite[Theorem~3.1]{Sar}. B\"or\"oczky and Kalantzopoulos \cite{BK} treated bodies sharing $n$ independent hyperplane symmetries. Van Handel \cite{vanHandel} proved that the local form of the logarithmic Brunn--Minkowski inequality holds for zonoids. Chen, Huang, Li, and Liu \cite{CHLL} proved the $L_p$-Brunn--Minkowski inequality for origin-symmetric convex bodies when $p<1$ is sufficiently close to $1$ and the equivalence of the local and global formulations for $0\leqslant p<1$. Putterman \cite{Putterman} gave another proof of this equivalence. Stancu \cite{StancuFlow} proposed a centro-affine flow approach to the logarithmic Minkowski inequality. For further developments, see \cite{BDStability,CLM,HenkPollehn,Iffland,KLlocal,KMBrascamp,KMSharp,KMBoundary,KM,Lu,Rotem,XiReverse,XiLeng} and the survey of B\"or\"oczky \cite{BoroczkySurvey}.

Very recently, an OpenAI preprint \cite{OpenAI} presented a proof of the logarithmic Brunn--Minkowski inequality for origin-symmetric convex bodies in all dimensions. For a discussion of the proof, see the lecture notes of Livshyts \cite{LivshytsCIRM}. Liu--Sun--Xiong \cite{LSX} extended the inequality to bodies with centroid at the origin and characterized equality.

\vskip 3pt
The other central topic of the Brunn-Minkowski theory, the Minkowski problem, is to characterize a geometric measure generated by convex bodies. The characterization of surface area measure is known as the classical Minkowski problem. From an analytic point of view, the surface area measure $S_K$ is the differential of volume under Minkowski addition: for $K,L\in\Kn$,
\[
 \frac{d}{dt}|K+tL|\Bigr|_{t=0^+} =\int_{\Sph}h_L\,dS_K.
\]
The study of the classical Minkowski problem has a long history and has led to a long series of highly influential works. See, e.g., \cite{AleksandrovMixedIII,CaffarelliLocalization,CaffarelliInterior,ChengYau,FenchelJessen,Lewy,MinkowskiVolume,Nirenberg,Pogorelov}.

The Minkowski problem for cone-volume measure is called the \emph{logarithmic Minkowski problem}. For $K\in\Ko$, its cone-volume measure $V_K$ is a finite Borel measure on $\Sph$, defined for Borel $\omega\subseteq\Sph$ by $$V_K(\omega)=\frac{1}{n}\int_\omega h_K(u)\,dS_K(u).$$
Geometrically, $V_K(\omega)$ is the volume swept out by the segments joining the origin to the boundary points of $K$ whose unique outer unit normal lies in $\omega$. From an analytic point of view, $V_K$ is the differential of volume under logarithmic Minkowski combination \cite{BLYZ12}: for $K,L\in\Ko$,
\[
 \frac{d}{dt} \bigl|\comb K L t\bigr|\Bigr|_{t=0^+} =n\int_{\Sph}\log\frac{h_L}{h_K}\,dV_K.
\]

In 2013, B\"or\"oczky--Lutwak--Yang--Zhang \cite{BLYZ13} introduced the celebrated \emph{subspace concentration condition} and proved that this condition is sufficient and necessary to guarantee the existence of solutions to the even logarithmic Minkowski problem. A nonzero finite Borel measure $\mu$ on $\Sph$ is said to satisfy this condition if
\[
 \frac{\mu(\Sph\cap \xi)}{\mu(\Sph)}\leqslant\frac{\dim \xi}{n}
\]
for each proper nonzero subspace $\xi$ of $\R^n$, and whenever equality holds for some $\xi$, there is a subspace $\xi'$ complementary to $\xi$ such that $\mu$ is concentrated on $\Sph\cap(\xi\cup \xi')$. For further results, see \cite{BHCentered,BHStability,HenkLinke,HuXiongPolytope,LSXDecomposition,LuSunXiongMixed,SunXiongMixed,SunXiongZonoids} on cone-volume measures and \cite{BHZ,CLZ,CLXTruncated,ChouWang,LLSX,StancuPlanar,SuXiongTruncated,ZhuLog} on the logarithmic Minkowski problem, as well as the surveys \cite{BoroczkySurvey,HYZSurvey}.

\vskip 3pt
Brunn--Minkowski inequalities extend beyond volume to other geometric functionals. Brascamp and Lieb \cite{BL} and Borell \cite{Borell} established the inequalities for the first Dirichlet eigenvalue and torsional rigidity, respectively. Colesanti \cite{Colesanti} characterized equality in both inequalities. For Newtonian capacity, Borell \cite{BorellCap} established the inequality, and Caffarelli, Jerison, and Lieb \cite{CJL} characterized equality; Colesanti and Salani \cite{CS} extended both results to $p$-capacity for $1<p<n$. Further examples include the Monge--Amp\`ere eigenvalue \cite{SalaniMA} and principal eigenvalues of a class of fully nonlinear homogeneous elliptic operators \cite{CFUnified}.

Minkowski problems associated with these functionals have also been studied. Jerison proved existence for the Minkowski problems associated with the first Dirichlet eigenvalue \cite[Theorem~7.4]{Jerison} and electrostatic capacity \cite{JerisonCap}. Colesanti et al.\ \cite{CNSXYZ} treated $p$-capacity, and Akman et al.\ \cite{AGHLV} considered a wider class of nonlinear capacities. Colesanti and Fimiani \cite{CF} solved the problem for torsional rigidity, with uniqueness up to translation following from the associated Brunn--Minkowski inequality and its equality case. Lutwak, Xi, Yang, and Zhang \cite[Theorems~1.1 and~1.2]{LXYZ} solved the existence part of the chord Minkowski problem for $q>0$. For $1\leqslant q\leqslant n+1$, they also gave a sufficient condition for an even measure to be the $q$th cone-chord measure of an origin-symmetric convex body.

In this paper, we study the logarithmic Brunn-Minkowski inequalities and the logarithmic Minkowski problems for four geometric functionals. For a convex body $K$, let $0<\lambda_1(K)\leqslant\lambda_2(K)\leqslant\cdots$ be the Dirichlet eigenvalues of $-\Delta$ on $\intt K$, counted with multiplicity. The first two functionals are the \emph{first Dirichlet eigenvalues} $\lambda_1(K)$ and the \emph{spectral zeta function} $$\zeta_\sigma(K)=\sum_{j=1}^\infty\lambda_j(K)^{-\sigma}$$ 
for $\sigma>\frac{n}{2}$. The third is the \emph{torsional rigidity} defined by $$\tau(K)=\int_Ku_K\,dx,$$ where $u_K$ solves $-\Delta u_K=1$ in $\intt K$ with zero boundary values. The last is the \emph{chord integral} defined for $q>-1$ by
\[
 I_q(K)=\int_{K\cap\ell\ne\varnothing}|K\cap\ell|^q\,d\ell,
\]
where $|K\cap\ell|$ is the length of the chord and $d\ell$ is the invariant measure on unoriented lines in $\R^n$.

The logarithmic variations of these four functionals give rise to the corresponding cone measures: for $K,L\in\Ko$,
\begin{align*}
 \frac{d}{dt}\lambda_1(\comb K L t)\Bigr|_{t=0^+} &=-\int_{\Sph}\log\frac{h_L}{h_K}\,d\lambda_1(K,u),\\
 \frac{d}{dt}\zeta_\sigma(\comb K L t)\Bigr|_{t=0^+} &=\int_{\Sph}\log\frac{h_L}{h_K}\,d\zeta_\sigma(K,u),\\
 \frac{d}{dt}\tau(\comb K L t)\Bigr|_{t=0^+} &=\int_{\Sph}\log\frac{h_L}{h_K}\,d\tau(K,u),\\
 \frac{d}{dt}I_q(\comb K L t)\Bigr|_{t=0^+} &=(n+q-1)\int_{\Sph}\log\frac{h_L}{h_K}\,dG_q(K,u).
\end{align*}
For $\lambda_1$ and $\tau$, these identities follow from the variation formulas of Jerison \cite[Theorem~7.5]{Jerison} and of Colesanti and Fimiani \cite[Corollary~4.2]{CF}, respectively. The formula for $I_q$ is a special case of the logarithmic variation formula of Lutwak, Xi, Yang, and Zhang \cite[(5.11)]{LXYZ}. We derive the formula for $\zeta_\sigma$ in Section~\ref{app:variations}.
\vskip 3pt
Our main results are as follows. We first prove the following general principle for geometric functionals. Write $K^m=K\times\cdots\times K$ for the $m$-fold Cartesian product of $K$. 

\begin{theorem}\label{thm:logbm-transfer}
Let $\beta\in\R$, $m_j\in\mathbb N_+$, and $\psi_j:\R^{nm_j}\to(-\infty,+\infty]$ be even convex functions for $j\in\mathbb N_+$. Suppose that for $K\in\Ke$ and $t>0$,
\[
H_K(t):=t^{\beta}\lim_{j\to\infty} \int_{(K/\sqrt{t})^{m_j}}e^{-\psi_j(x)}\,dx\in (0,\infty).
\]
\begin{enumerate}
\item[(i)] Let $\alpha\in\R$. If $\mathcal J_\alpha(K):=\int_0^\infty t^{\alpha-1}H_K(t)\,dt<\infty$ for every $K\in\Ke$, then for $K,L\in\Ke$ and $0<\theta<1$,
\[
 \mathcal J_\alpha\bigl(\comb K L\theta\bigr)\geqslant \mathcal J_\alpha(K)^{1-\theta}\mathcal J_\alpha(L)^\theta,
\]
with equality if and only if there exist $r,s>0$ such that $H_L(rt)=sH_K(t)$ and $H_{\comb K L\theta}(r^\theta t)=s^\theta H_K(t)$ for $t>0$.
\item[(ii)] If $\Lambda(K):=-\lim\limits_{t\to\infty}t^{-1}\log H_K(t)\in(0,\infty)$ for every $K\in\Ke$, then for $K,L\in\Ke$ and $0<\theta<1$,
\[
 \Lambda\bigl(\comb K L\theta\bigr)\leqslant \Lambda(K)^{1-\theta}\Lambda(L)^\theta.
\]
\end{enumerate}
\end{theorem}

Applying Theorem~\ref{thm:logbm-transfer} to different $H_K$, we prove the four logarithmic Brunn--Minkowski inequalities below. Specifically,
\begin{alignat*}{2}
 H_K(t)&=\sum_{j=1}^\infty e^{-t\lambda_j(K)} &\quad&\Longrightarrow\quad \Lambda(K)=\lambda_1(K),\  \mathcal J_\sigma(K)=\Gamma(\sigma)\zeta_\sigma(K),\\
 H_K(t)&=\int_{K^2}k_t^K(x,y)\,dx\,dy &&\Longrightarrow\quad \mathcal J_1(K)=\tau(K),\\
 H_K(t)&=\int_{K^2}e^{-\frac{|x-y|^2}{4t}}\,dx\,dy &&\Longrightarrow\quad \mathcal J_{\frac{q-n-1}{2}}(K)=c_{n,q}I_q(K).
\end{alignat*}
Here $k_t^K$ is the Dirichlet heat kernel on $\intt K$, and $c_{n,q}>0$ depends only on $n$ and $q$.

\begin{theorem}
\label{thm:geometric-logbm}
Let $K,L\in\Ke$, $0<\theta<1$, $\sigma>\frac{n}{2}$, and $1<q<n+1$. Then
\begin{align*}
 \lambda_1(\comb K L\theta)&\leqslant \lambda_1(K)^{1-\theta}\lambda_1(L)^\theta,\\
 \zeta_\sigma(\comb K L\theta)&\geqslant \zeta_\sigma(K)^{1-\theta}\zeta_\sigma(L)^\theta,\\
 \tau(\comb K L\theta)&\geqslant \tau(K)^{1-\theta}\tau(L)^\theta,\\
 I_q(\comb K L\theta)&\geqslant I_q(K)^{1-\theta}I_q(L)^\theta.
\end{align*}
In each case, equality holds if and only if $K$ and $L$ are dilates.
\end{theorem}

We also solve the corresponding logarithmic Minkowski problems for even measures.

\begin{theorem}
\label{thm:geometric-inverse}
Let $\mu$ be a nonzero finite even Borel measure on $\Sph$. In the inequalities below, $\xi$ ranges over all nonzero proper subspaces of $\R^n$.
\begin{enumerate}
\renewcommand{\labelenumi}{(\roman{enumi})}
\item There exists $K\in\Ke$ such that $\lambda_1(K,\cdot)=\mu$ if and only if
\begin{equation*}
 \frac{\mu(\Sph\cap \xi)}{\mu(\Sph)} <1.
\end{equation*}
\item For $\sigma>\frac{n}{2}$, there exists $K\in\Ke$ such that $\zeta_\sigma(K,\cdot)=\mu$ if and only if
\begin{equation*}
 \frac{\mu(\Sph\cap \xi)}{\mu(\Sph)} <1-\frac{n-\dim\xi}{2\sigma}.
\end{equation*}
\item There exists $K\in\Ke$ such that $\tau(K,\cdot)=\mu$ if and only if
\begin{equation*}
 \frac{\mu(\Sph\cap \xi)}{\mu(\Sph)} <\frac{\dim\xi+2}{n+2}.
\end{equation*}
\item For $1<q<n+1$, there exists $K\in\Ke$ such that $G_q(K,\cdot)=\mu$ if and only if
\begin{equation*}
 \frac{\mu(\Sph\cap \xi)}{\mu(\Sph)} <\frac{\dim\xi+\min\{\dim\xi,q-1\}}{n+q-1}.
\end{equation*}
\end{enumerate}
Moreover, the solution is unique in $\Ke$.
\end{theorem}

For origin-symmetric convex bodies, Theorem~\ref{thm:geometric-logbm} strengthens the corresponding Brunn--Minkowski inequalities, since $\comb K L\theta\subseteq(1-\theta)K+\theta L$ and the four functionals are monotone under inclusion and homogeneous.

Letting $q$ tend to $1$ or $n+1$ in the chord integral inequality of Theorem~\ref{thm:geometric-logbm} recovers the logarithmic Brunn--Minkowski inequality for volume. Indeed, $I_q(K)$ is continuous in $q$, with $I_1(K)=|K|$ and $I_{n+1}(K)=\frac{n+1}{\omega_n}|K|^2$, where $\omega_n$ is the volume of the Euclidean unit ball.

Theorems~\ref{thm:geometric-logbm} and~\ref{thm:geometric-inverse}~(iii) answer two questions of Crasta and Fragal\`a \cite[(11)--(12)]{CrF}, which state that on $\Ke$, torsional rigidity satisfies the logarithmic Brunn--Minkowski inequality and the cone torsion measure uniquely determines the body.

The inequalities in Theorem~\ref{thm:geometric-inverse} are sharp for each fixed $\xi$ (see Remark~\ref{rem:sharpness-existence}). For the first Dirichlet eigenvalue and torsional rigidity, Langharst and Ulivelli \cite{LU} proved existence for even measures satisfying the stronger condition $\frac{\mu(\Sph\cap\xi)}{\mu(\Sph)}<\frac{\dim\xi}{n}$. For chord integrals, Lutwak, Xi, Yang, and Zhang \cite{LXYZ} proved that the condition in (iv) is sufficient for even measures. Guo, Xi, and Zhao \cite{GXZ} removed the evenness assumption. For $q\in\{2,\ldots,n\}$, Lutwak, Xi, Yang, and Zhang obtained a necessary upper bound and proved it sharp when $\dim\xi<q-1$. Part~(iv) shows that for $1<q<n+1$ the condition is also necessary for an origin-symmetric solution and that this solution is unique.

The paper is organized as follows. In Section~\ref{sec:preliminaries}, we collect some basic facts about convex bodies, the four geometric functionals, and heat kernels. We prove the inequalities in Theorems~\ref{thm:logbm-transfer} and~\ref{thm:geometric-logbm} in Section~\ref{sec:transfer} and characterize equality in Section~\ref{sec:rigidity-common}. We establish a variational principle for logarithmic Minkowski problems in Section~\ref{sec:abstract-inverse} and apply it in Section~\ref{sec:inverse-applications} to prove Theorem~\ref{thm:geometric-inverse}.

\vskip 25pt
\section{\bf Preliminaries}\label{sec:preliminaries}
\vskip 5pt

\subsection{Convex bodies and Wulff shapes}\label{sec:convex-wulff}

For background, see Gardner \cite{Gardner}, Hug--Weil \cite{HugWeil}, and Schneider \cite{Schneider}.

As usual, $x\cdot y$ denotes the standard inner product of $x$ and $y$, and $|x|$ denotes the standard Euclidean norm of $x$ in $\R^n$. Write $\Sph=\{x\in\R^n:|x|=1\}$ for the boundary of the Euclidean unit ball $B^n$ in $\R^n$. The volume of $B^n$ is $\omega_n=\frac{\pi^{\frac{n}{2}}}{\Gamma(1+\frac{n}{2})}$.

For a linear subspace $\xi$ of $\R^n$, let $\xi^\perp$ denote its orthogonal complement and $\mathrm{P}_\xi:\R^n\to\xi$ the orthogonal projection onto $\xi$. Write $B_\xi=\xi\cap B^n$.

For $K\in\Kn$, its \emph{support function} $h_K:\R^n\to\R$ is defined by
\[
 h_K(x)=\max\{x\cdot y:y\in K\},\qquad x\in\R^n.
\]

$\Kn$ is equipped with the \emph{Hausdorff metric} $\delta_H$, defined by
\[
 \delta_H(K,L)=\max_{u\in\Sph}|h_K(u)-h_L(u)|,\qquad K,L\in\Kn.
\]

The \emph{surface area measure} $S_K$ of $K\in\Kn$ is a finite Borel measure on $\Sph$, defined for Borel sets $\omega\subseteq\Sph$ by $S_K(\omega)=\mathcal H^{n-1}(\nu_K^{-1}(\omega))$, where $\mathcal H^{n-1}$ is the $(n-1)$-dimensional Hausdorff measure, and $\nu_K:\partial'K\to\Sph$ is the Gauss map of $K$, defined on $\partial'K$, the set of points of $\partial K$ that have a unique outer unit normal. Recall that $\mathcal H^{n-1}(\partial K\setminus\partial'K)=0$. See, e.g., \cite[p.~84]{Schneider} for details.

Denote by $C(\Sph)$ the space of continuous real-valued functions on $\Sph$, and by $C_e(\Sph)$ its subspace of even functions.

For positive $f\in C(\Sph)$, the Wulff shape (or Aleksandrov body) determined by $f$ is
\[
 [f]=\{x\in\R^n:x\cdot u\leqslant f(u)\text{ for }u\in\Sph\}.
\]
Then $h_{[f]}\leqslant f$ on $\Sph$, and $[h_K]=K$ for $K\in\Ko$. 

\begin{lemma}\label{lem:wulff-comparison}
If $\mathcal F:\Ko\to(0,\infty)$ is monotone with respect to inclusion and homogeneous of degree one, then for positive $f,g\in C(\Sph)$,
\begin{equation}\label{eq:wulff-ratio-estimate}
 \min_{\Sph}\frac{f}{g}\leqslant\frac{\mathcal F([f])}{\mathcal F([g])}\leqslant\max_{\Sph}\frac{f}{g}.
\end{equation}
In particular, $\mathcal F$ is continuous on $\Ko$. If $\mathcal F:\Ke\to(0,\infty)$ satisfies the same assumptions, then \eqref{eq:wulff-ratio-estimate} holds for positive $f,g\in C_e(\Sph)$, and $\mathcal F$ is continuous on $\Ke$.
\end{lemma}
\begin{proof}
Write $m=\min_{\Sph}\frac{f}{g}$ and $M=\max_{\Sph}\frac{f}{g}$. Then $mg\leqslant f\leqslant Mg$. By the definition of the Wulff shape, we have $m[g]\subseteq[f]\subseteq M[g]$. Monotonicity and homogeneity of $\mathcal F$ then yield \eqref{eq:wulff-ratio-estimate}.

If $K_j\to K$ in $\Ko$, then $\frac{h_{K_j}}{h_K}\to 1$ uniformly on $\Sph$. Taking $f=h_{K_j}$ and $g=h_K$ in \eqref{eq:wulff-ratio-estimate}, we obtain $\mathcal F(K_j)\to\mathcal F(K)$.

For even $f$ and $g$, the Wulff shapes $[f]$ and $[g]$ are origin-symmetric, so the same proof applies on $\Ke$.
\end{proof}

For $K,L\in\Ko$ and $0<\theta<1$, their logarithmic Minkowski combination is defined by $\comb K L\theta=[h_K^{1-\theta}h_L^\theta]$. By the arithmetic--geometric mean inequality,
\[
 \comb K L\theta\subseteq[(1-\theta)h_K+\theta h_L]=(1-\theta)K+\theta L.
\]
For $m\in\mathbb N_+$, let $K^m=K\times\cdots\times K$ denote the $m$-fold Cartesian product.

\begin{lemma}\label{lem:product-logarithmic-combination}
For $K,L\in\Ko$, $m\in\mathbb N_+$, and $0<\theta<1$,
\[
 \comb{(K^m)}{(L^m)}\theta =(\comb K L\theta)^m.
\]
\end{lemma}
\begin{proof}
First, observe that the support function of $K^m$ is $h_{K^m}(u_1,\ldots,u_m)=\sum_{i=1}^m h_K(u_i)$.

Assume $x=(x_1,\ldots,x_m)\in\comb{(K^m)}{(L^m)}\theta$. For $i=1,\ldots,m$ and $v\in\Sph$, take $u=(0,\ldots,0,v,0,\ldots,0)$, with $v$ in the $i$th position. By the definition of the logarithmic Minkowski combination,
\[
 x_i\cdot v=x\cdot u\leqslant h_{K^m}(u)^{1-\theta}h_{L^m}(u)^\theta=h_K(v)^{1-\theta}h_L(v)^\theta.
\]
Thus $x\in(\comb K L\theta)^m$. Hence $\comb{(K^m)}{(L^m)}\theta\subseteq(\comb K L\theta)^m$.

Conversely, assume $x=(x_1,\ldots,x_m)\in(\comb K L\theta)^m$. By definition and the positive homogeneity of support functions, $x_i\cdot v\leqslant h_K(v)^{1-\theta}h_L(v)^\theta$ for $v\in\R^n$ and $i=1,\ldots,m$. Hence, for $u=(u_1,\ldots,u_m)\in\mathbb S^{nm-1}$, H\"older's inequality yields
\begin{align*}
 x\cdot u&=\sum_{i=1}^m x_i\cdot u_i\leqslant\sum_{i=1}^m h_K(u_i)^{1-\theta}h_L(u_i)^\theta\\
 &\leqslant\Bigl(\sum_{i=1}^m h_K(u_i)\Bigr)^{1-\theta}\Bigl(\sum_{i=1}^m h_L(u_i)\Bigr)^\theta=h_{K^m}(u)^{1-\theta}h_{L^m}(u)^\theta.
\end{align*}
Thus $x\in\comb{(K^m)}{(L^m)}\theta$. Hence $(\comb K L\theta)^m\subseteq\comb{(K^m)}{(L^m)}\theta$.
\end{proof}

By Saroglou \cite[Theorem~3.1]{Sar}, the logarithmic Brunn--Minkowski inequality for volume implies the following inequality for even log-concave functions.

\begin{lemma}\label{lem:saroglou-density}
If $\psi:\R^n\to(-\infty,+\infty]$ is even and convex, then for $K,L\in\Ke$ and $0<\theta<1$,
\[
 \int_{\comb K L\theta}e^{-\psi(x)}\,dx \geqslant \Bigl(\int_K e^{-\psi(x)}\,dx\Bigr)^{1-\theta} \Bigl(\int_L e^{-\psi(x)}\,dx\Bigr)^\theta.
\]
\end{lemma}

Note that the integrals in Lemma~\ref{lem:saroglou-density} are finite, since evenness and convexity imply that $e^{-\psi(x)}\leqslant e^{-\psi(0)}<\infty$.

We also need the following one-dimensional Pr\'ekopa--Leindler inequality \cite{Leindler1972,Prekopa1971}.

\begin{lemma}\label{lem:prekopa-leindler}
Let $f,g,h$ be nonnegative integrable functions on $\R$ and $0<\theta<1$. If
\[
 h((1-\theta)x+\theta y)\geqslant f(x)^{1-\theta}g(y)^\theta
\]
for $x,y\in\R$, then
\[
 \int_\R h(x)\,dx\geqslant \Bigl(\int_\R f(x)\,dx\Bigr)^{1-\theta} \Bigl(\int_\R g(x)\,dx\Bigr)^\theta.
\]
\end{lemma}

\subsection{Geometric functionals}\label{sec:geometric-functionals}

For $K\in\Kn$, let $0<\lambda_1(K)\leqslant\lambda_2(K)\leqslant\cdots$ denote the eigenvalues of $-\Delta$ on $\intt K$ with zero Dirichlet boundary values, counted with multiplicity. For $K,L\in\Kn$ with $K\subseteq L$, $j\geqslant1$, and $r>0$,
\begin{equation}\label{eq:eigenvalue-monotonicity-scaling}
 \lambda_j(K)\geqslant\lambda_j(L),\qquad \lambda_j(rK)=r^{-2}\lambda_j(K).
\end{equation}
See, e.g., \cite[Section~8.12]{GilbargTrudinger}.

For $\sigma>\frac{n}{2}$, the spectral zeta function is $\zeta_\sigma(K)=\sum_{j=1}^\infty\lambda_j(K)^{-\sigma}$. If $K\subseteq L$, then $\zeta_\sigma(K)\leqslant\zeta_\sigma(L)$ by \eqref{eq:eigenvalue-monotonicity-scaling}. The series converges by \eqref{eq:zeta-heat-representation}.

The torsional rigidity is $\tau(K)=\int_Ku_K\,dx$, where $u_K\in H_0^1(\intt K)$ is the unique weak solution of $-\Delta u_K=1$ in $\intt K$. Its variational characterization implies $\tau(K)\leqslant\tau(L)$ for $K\subseteq L$. Since $u_{rK}(rx)=r^2u_K(x)$ for $r>0$, we have $\tau(rK)=r^{n+2}\tau(K)$. See, e.g., \cite[Section~2]{CF}.

Let $\mathscr{L}^n$ denote the affine Grassmannian of unoriented lines in $\R^n$. For $K\in\Kn$ and $q>-1$, the chord integral is
\[
 I_q(K)=\int_{K\cap\ell\ne\varnothing}|K\cap\ell|^q\,d\ell,
\]
where $|K\cap\ell|$ denotes the length of the chord and $d\ell$ is the invariant measure normalized so that $I_1(B^n)=\omega_n$. For almost every line $\ell$ meeting $K$, $|K\cap\ell|>0$. The integral $I_q(K)$ is finite for $q>-1$ \cite[Lemma~3.1]{LXYZ}.

If $K\subseteq L$, then $I_q(K)\leqslant I_q(L)$ for $q>0$. For $q\in[a,b]\subseteq(-1,\infty)$ and $r=|K\cap\ell|>0$, we have $r^q\leqslant r^a+r^b$. Since $I_a(K)+I_b(K)<\infty$, dominated convergence shows that $q\mapsto I_q(K)$ is continuous on $(-1,\infty)$ for fixed $K$.

The Crofton formula and the Poincar\'e--Hadwiger formula \cite{Ren,Santalo} give
\[
 I_1(K)=|K|,\qquad I_{n+1}(K)=\frac{n+1}{\omega_n}|K|^2.
\]

For $1<q<n+1$, the chord integral has the classical representation
\begin{equation}\label{ri:eq:chord}
 I_q(K)=\frac{q(q-1)}{n\omega_n}E_{n+1-q}(K),
\end{equation}
where $E_a(K)=\int_{K^2}|x-y|^{-a}\,dx\,dy$ is the Riesz energy of $K$ for $0<a<n$. See \cite{Ren} and \cite[p.~238, (14.25)]{Santalo}.

For $\sigma>\frac{n}{2}$ and $1<q<n+1$, the normalized functionals $\lambda_1^{-\frac{1}{2}}$, $\zeta_\sigma^{\frac{1}{2\sigma}}$, $\tau^{\frac{1}{n+2}}$, and $I_q^{\frac{1}{n+q-1}}$ are positive, monotone with respect to inclusion, and homogeneous of degree one. Their continuity on $\Ko$ follows from Lemma~\ref{lem:wulff-comparison}.

\subsection{Heat kernels}\label{sec:heat-kernels}

The Euclidean heat kernel is
\begin{equation}\label{eq:euclidean-heat-kernel}
 p_t(x,y)=(4\pi t)^{-\frac{n}{2}} \exp\!\Bigl(-\frac{|x-y|^2}{4t}\Bigr), \qquad t>0,\ x,y\in\R^n.
\end{equation}
For $K\in\Kn$, let $k_t^K$, $t>0$, denote the Dirichlet heat kernel on $\intt K$ for the equation $\partial_tu=\Delta u$. That is, for $f\in L^2(K)$, the function
\[
 u(t,x)=\int_Kk_t^K(x,y)f(y)\,dy, \qquad t>0,x\in\intt K,
\]
solves this equation with zero Dirichlet boundary values, and $u(t,\cdot)\to f$ in $L^2(K)$ as $t\to0^+$. See \cite[Chapters~4 and 6]{Ouhabaz} for details.

Let $\{\phi_j\}_{j\geqslant1}$ be a real orthonormal basis of $L^2(K)$ consisting of Dirichlet eigenfunctions. Each $\phi_j\in H_0^1(\intt K)$ satisfies $-\Delta\phi_j=\lambda_j(K)\phi_j$ in the weak sense. Since $\intt K$ is connected, $\lambda_1(K)$ is simple, and we choose $\phi_1$ strictly positive in $\intt K$. For these spectral facts, see \cite[Theorems~8.37 and 8.38]{GilbargTrudinger}. The kernel $k_t^K$ has the spectral expansion
\begin{equation}\label{eq:heat-kernel-spectral}
 k_t^K(x,y)=\sum_{j=1}^\infty e^{-t\lambda_j(K)}\phi_j(x)\phi_j(y), \quad \text{for }t>0,\ x,y\in\intt K.
\end{equation}
See \cite[Lemma~2.1]{DaviesSimon}. In particular, $k_t^K(x,y)=k_t^K(y,x)$.

For $t>0$, $x,y\in\intt K$, and integers $j\geqslant2$, set $x_0=x$, $x_j=y$, and define the approximating kernels
\begin{equation}\label{eq:heat-kernel-product}
 k_{t,j}^K(x,y)=\int_{K^{j-1}} \prod_{i=1}^{j}p_{\frac{t}{j}}(x_{i-1},x_i)\,dx_1\cdots dx_{j-1}.
\end{equation}

The following approximation result is proved in \cite[p.~384, (6.4)]{BL}.

\begin{lemma}\label{lem:heat-kernel-product}
Let $K\in\Kn$, $t>0$ and $x,y\in\intt K$. Then $0\leqslant k_{t,j}^K(x,y)\leqslant p_t(x,y)$, $j\geqslant2$, and $k_{t,j}^K(x,y)\longrightarrow k_t^K(x,y)$ as $j\to\infty$.
\end{lemma}

If $K\subseteq L$, then \eqref{eq:heat-kernel-product} gives $k_{t,j}^K(x,y)\leqslant k_{t,j}^L(x,y)$ for $t>0$ and $x,y\in\intt K$. Letting $j\to\infty$, we obtain $k_t^K(x,y)\leqslant k_t^L(x,y)$.

\vskip 25pt
\section{\bf Logarithmic Brunn--Minkowski inequalities for geometric functionals}\label{sec:transfer}
\vskip 5pt

In this section, we prove the inequalities in Theorems~\ref{thm:logbm-transfer} and~\ref{thm:geometric-logbm}. The equality conditions are established in Section~\ref{sec:rigidity-common}.

\subsection{A transfer principle}

\begin{theorem}\label{thm:transfer-inequalities}
Let $\beta\in\R$, $m_j\in\mathbb N_+$, and $\psi_j:\R^{nm_j}\to(-\infty,+\infty]$ be even convex functions for $j\in\mathbb N_+$. Suppose that for $K\in\Ke$ and $t>0$,
\begin{equation}\label{eq:transfer-profile-representation}
 H_K(t):=t^{\beta}\lim_{j\to\infty} \int_{(K/\sqrt{t})^{m_j}}e^{-\psi_j(x)}\,dx\in (0,\infty).
\end{equation}
\begin{enumerate}
\item[(i)] Let $\alpha\in\R$. If $\mathcal J_\alpha(K):=\int_0^\infty t^{\alpha-1}H_K(t)\,dt<\infty$ for every $K\in\Ke$, then for $K,L\in\Ke$ and $0<\theta<1$,
\begin{equation}\label{eq:transfer-integral-logbm}
 \mathcal J_\alpha\bigl(\comb K L\theta\bigr)\geqslant \mathcal J_\alpha(K)^{1-\theta}\mathcal J_\alpha(L)^\theta.
\end{equation}
\item[(ii)] If $\Lambda(K):=-\lim\limits_{t\to\infty}t^{-1}\log H_K(t)\in(0,\infty)$ for every $K\in\Ke$, then for $K,L\in\Ke$ and $0<\theta<1$,
\[
 \Lambda\bigl(\comb K L\theta\bigr)\leqslant \Lambda(K)^{1-\theta}\Lambda(L)^\theta.
\]
\end{enumerate}
\end{theorem}
\begin{proof}
Write $M=\comb K L\theta$ for $K,L\in\Ke$ and $0<\theta<1$. We first prove that for $t_0,t_1>0$,
\begin{equation}\label{eq:abstract-joint-time}
 H_M(t_0^{1-\theta}t_1^\theta) \geqslant H_K(t_0)^{1-\theta}H_L(t_1)^\theta.
\end{equation}

To this end, let $m\in\mathbb N_+$ and let $\psi:\R^{nm}\to(-\infty,+\infty]$ be an even convex function. Applying Lemma~\ref{lem:saroglou-density} in dimension $nm$ to $K^m,L^m$ and using Lemma~\ref{lem:product-logarithmic-combination} gives
\begin{equation}\label{eq:weighted-product-comparison}
 \int_{M^m}e^{-\psi(x)}\,dx \geqslant \Bigl(\int_{K^m}e^{-\psi(x)}\,dx\Bigr)^{1-\theta} \Bigl(\int_{L^m}e^{-\psi(x)}\,dx\Bigr)^\theta.
\end{equation}

By definition,
\[
 \comb{(K/\sqrt{t_0})}{(L/\sqrt{t_1})}\theta=M/\sqrt{t_0^{1-\theta}t_1^\theta}.
\]
Applying \eqref{eq:weighted-product-comparison} to these two bodies with $m=m_j$ and $\psi=\psi_j$, Letting $j\to\infty$ and using \eqref{eq:transfer-profile-representation}, we conclude that
\[
 \frac{H_M(t_0^{1-\theta}t_1^\theta)} {(t_0^{1-\theta}t_1^\theta)^{\beta}} \geqslant \Bigl(\frac{H_K(t_0)}{t_0^{\beta}}\Bigr)^{1-\theta} \Bigl(\frac{H_L(t_1)}{t_1^{\beta}}\Bigr)^\theta.
\]

For part~(i), write $f_D(s)=e^{\alpha s}H_D(e^s)$ for $D=K,L,M$. Then
\[
 \int_\R f_D(s)\,ds=\int_\R e^{\alpha s}H_D(e^s)\,ds=\int_0^\infty t^{\alpha-1}H_D(t)\,dt=\mathcal J_\alpha(D)<\infty.
\]
By \eqref{eq:abstract-joint-time}, we have
\[
 f_M((1-\theta)s_0+\theta s_1) \geqslant f_K(s_0)^{1-\theta}f_L(s_1)^\theta
\]
for $s_0,s_1\in\R$. Lemma~\ref{lem:prekopa-leindler} now implies the inequality in part~(i).

For part~(ii), fix $t>0$ and take $t_0=\frac{t}{\Lambda(K)}$ and $t_1=\frac{t}{\Lambda(L)}$ in \eqref{eq:abstract-joint-time}. Then
\[
 -\frac1t\log H_M\Bigl( \frac{t}{\Lambda(K)^{1-\theta}\Lambda(L)^\theta} \Bigr) \leqslant -(1-\theta)\frac{\log H_K(t_0)}t -\theta\frac{\log H_L(t_1)}t.
\]
Letting $t\to\infty$ and using the definition of $\Lambda$ gives
\[
 \frac{\Lambda(M)} {\Lambda(K)^{1-\theta}\Lambda(L)^\theta} \leqslant (1-\theta)\cdot1+\theta\cdot1=1,
\]
as desired.
\end{proof}

\subsection{Inequalities for geometric functionals}\label{sec:functional-inequalities}

\begin{theorem}\label{thm:functional-inequalities}
Let $K,L\in\Ke$, $0<\theta<1$, $\sigma>\frac{n}{2}$, and $1<q<n+1$. Then
\begin{align}
 \lambda_1(\comb K L\theta)&\leqslant \lambda_1(K)^{1-\theta}\lambda_1(L)^\theta,\label{eq:target-eigen-logbm}\\
 \zeta_\sigma(\comb K L\theta)&\geqslant \zeta_\sigma(K)^{1-\theta}\zeta_\sigma(L)^\theta,\label{eq:target-zeta-logbm}\\
 \tau(\comb K L\theta)&\geqslant \tau(K)^{1-\theta}\tau(L)^\theta,\label{eq:target-torsion-logbm}\\
 I_q(\comb K L\theta)&\geqslant I_q(K)^{1-\theta}I_q(L)^\theta.\label{eq:target-chord-logbm}
\end{align}
\end{theorem}
\begin{proof}
Let $k_t^K$ and $\{\phi_j\}_{j\geqslant1}$ be the Dirichlet heat kernel and eigenbasis in Section~\ref{sec:heat-kernels}. For $t>0$, define
\begin{align}
 Q_K(t)&=\int_{K^2}k_t^K(x,y)\,dx\,dy,\label{eq:heat-content-definition}\\
 Z_K(t)&=\int_Kk_t^K(x,x)\,dx,\label{eq:heat-trace-definition}\\
 F_K(t)&=\int_{K^2}\exp\!\Bigl(-\frac{|x-y|^2}{4t}\Bigr)\,dx\,dy.\label{eq:gaussian-integral-definition}
\end{align}

We apply Theorem~\ref{thm:transfer-inequalities} with the choices below, using part~(ii) for $\lambda_1$ and part~(i) for the other three functionals.
\begin{center}
\renewcommand{\arraystretch}{1.25}
\setlength{\tabcolsep}{9pt}
\begin{tabular}{lc@{\hspace{24pt}}c@{\hspace{24pt}}c}
\toprule
Functional & $H_K$ & $\beta$ & $\alpha$ \\
\midrule
First Dirichlet eigenvalue $\lambda_1$ & $Z_K$ & $0$ & --- \\
Spectral zeta function $\zeta_\sigma$ & $Z_K$ & $0$ & $\sigma$ \\
Torsional rigidity $\tau$ & $Q_K$ & $\dfrac{n}{2}$ & $1$ \\
Chord integral $I_q$ & $F_K$ & $n$ & $\dfrac{q-n-1}{2}$ \\
\bottomrule
\end{tabular}
\end{center}

We first express $Q_K$ and $Z_K$ in the form \eqref{eq:transfer-profile-representation}.

Fix $t>0$ and let $k_{t,j}^K$ be given by \eqref{eq:heat-kernel-product}. For $x,y\in\intt K$, Lemma~\ref{lem:heat-kernel-product} gives
\begin{equation}\label{eq:heat-kernel-limit-bound}
\begin{aligned}
 k_{t,j}^K(x,y)&\longrightarrow k_t^K(x,y)\qquad\text{as }j\to\infty,\\
 0\leqslant k_{t,j}^K(x,y)&\leqslant p_t(x,y)\leqslant(4\pi t)^{-\frac{n}{2}},\; j\geqslant2.
\end{aligned}
\end{equation}
Since $|\partial K|=0$, it follows from the dominated convergence theorem, \eqref{eq:heat-content-definition} and \eqref{eq:heat-trace-definition} that
\begin{align*}
 Q_K(t)&=\lim_{j\to\infty}\int_{K^2}k_{t,j}^K(x,y)\,dx\,dy,\\
 Z_K(t)&=\lim_{j\to\infty}\int_Kk_{t,j}^K(x,x)\,dx.
\end{align*}

Since $p_{\frac{t}{j}}(\sqrt t\,x,\sqrt t\,y)=t^{-\frac{n}{2}}p_{\frac{1}{j}}(x,y)$ by \eqref{eq:euclidean-heat-kernel}, substituting \eqref{eq:heat-kernel-product} into these integrals and changing variables by $x_i=\sqrt t\,y_i$ gives
\begin{align*}
 Q_K(t)&=t^{\frac{n}{2}}\lim_{j\to\infty} \int_{(K/\sqrt{t})^{j+1}}\prod_{i=1}^j p_{\frac{1}{j}}(y_{i-1},y_i)\,dy_0\cdots dy_j,\\
 Z_K(t)&=\lim_{j\to\infty} \int_{(K/\sqrt{t})^j}\prod_{i=1}^j p_{\frac{1}{j}}(y_{i-1},y_i)\,dy_0\cdots dy_{j-1},
\end{align*}
where $y_j=y_0$ for $Z_K$. By \eqref{eq:euclidean-heat-kernel}, the negative logarithm of each integrand is
\[
 \psi_j(y_0,\ldots,y_j)=\frac j4\sum_{i=1}^j|y_i-y_{i-1}|^2+\frac{nj}{2}\log\frac{4\pi}{j}.
\]
Since $\psi_j$ is even, convex, and independent of $K,t$, it follows that \eqref{eq:transfer-profile-representation} holds with $m_j=j+1$ and $\beta=\frac{n}{2}$ for $Q_K$, and with $m_j=j$ and $\beta=0$ for $Z_K$.
\vskip 3pt
\textbf{Case 1.} Prove the inequality \eqref{eq:target-eigen-logbm} for the first Dirichlet eigenvalue.

By Step~1 and part~(ii) of Theorem~\ref{thm:transfer-inequalities}, it suffices to prove that
\[
 \lambda_1(K)=-\lim_{t\to\infty}\frac{\log H_K(t)}{t}
\]
with either $H_K=Z_K$ or $H_K=Q_K$.

For fixed $t>0$, the summands in \eqref{eq:heat-kernel-spectral} with $y=x$ satisfy $e^{-t\lambda_j(K)}\phi_j(x)^2\geqslant0$. Thus, applying monotone convergence in \eqref{eq:heat-trace-definition} and using $\int_K\phi_j^2\,dx=1$, we obtain
\begin{equation}\label{eq:heat-trace-expansion}
 Z_K(t)=\int_K k_t^K(x,x)\,dx=\sum_{j=1}^\infty e^{-t\lambda_j(K)}\int_K\phi_j(x)^2\,dx=\sum_{j=1}^\infty e^{-t\lambda_j(K)}.
\end{equation}
Letting $j\to\infty$ in \eqref{eq:heat-kernel-limit-bound} yields $k_t^K(x,x)\leqslant p_t(x,x)=(4\pi t)^{-\frac{n}{2}}$. Hence
\begin{equation}\label{eq:heat-trace-upper-bound}
 Z_K(t)\leqslant\int_Kp_t(x,x)\,dx=(4\pi t)^{-\frac{n}{2}}|K|<\infty.
\end{equation}

By the Cauchy--Schwarz inequality, we have $\int_K|\phi_j|\,dx\leqslant |K|^{\frac{1}{2}}$. Thus
\[
 \sum_{j=1}^\infty e^{-t\lambda_j(K)}\int_{K^2}|\phi_j(x)\phi_j(y)|\,dx\,dy=\sum_{j=1}^\infty e^{-t\lambda_j(K)}\Bigl(\int_K|\phi_j|\,dx\Bigr)^2\leqslant |K|Z_K(t)<\infty.
\]
This, together with Fubini's theorem, \eqref{eq:heat-kernel-spectral}, and \eqref{eq:heat-content-definition}, yields
\begin{equation}\label{eq:heat-content-expansion}
 Q_K(t)=\int_{K^2}k_t^K(x,y)\,dx\,dy=\sum_{j=1}^\infty e^{-t\lambda_j(K)}\Bigl(\int_K\phi_j(x)\,dx\Bigr)^2.
\end{equation}

For $t\geqslant1$, the inequality $\lambda_j(K)\geqslant\lambda_1(K)$ implies
\[
 e^{-t\lambda_j(K)} =e^{-\lambda_j(K)}e^{-(t-1)\lambda_j(K)} \leqslant e^{-\lambda_j(K)}e^{-(t-1)\lambda_1(K)}.
\]
This, together with \eqref{eq:heat-trace-expansion} and \eqref{eq:heat-content-expansion}, yields
\begin{align}
 Z_K(t)&\leqslant Z_K(1)e^{-(t-1)\lambda_1(K)},\label{eq:heat-trace-exponential-bound}\\
 Q_K(t)&\leqslant Q_K(1)e^{-(t-1)\lambda_1(K)}.\label{eq:heat-content-exponential-bound}
\end{align}

Retaining only the first term in \eqref{eq:heat-trace-expansion} and \eqref{eq:heat-content-expansion} gives
\[
 Z_K(t)\geqslant e^{-t\lambda_1(K)}, \qquad Q_K(t)\geqslant e^{-t\lambda_1(K)}\Bigl(\int_K\phi_1(x)\,dx\Bigr)^2.
\]
Since $\phi_1>0$ in $\intt K$ (see Section~\ref{sec:heat-kernels}), these lower bounds and \eqref{eq:heat-trace-exponential-bound}--\eqref{eq:heat-content-exponential-bound} imply
\begin{equation}\label{eq:heat-profile-decay}
 -\lim_{t\to\infty}\frac{\log Z_K(t)}t =-\lim_{t\to\infty}\frac{\log Q_K(t)}t=\lambda_1(K).
\end{equation}
Inequality~\eqref{eq:target-eigen-logbm} now follows from Theorem~\ref{thm:transfer-inequalities}~(ii).
\vskip 3pt
\textbf{Case 2.} Prove the inequality \eqref{eq:target-zeta-logbm} for spectral zeta functions.

Fix $j\geqslant1$. Since $\lambda_j(K)>0$, it follows that
\begin{align*}
 \Gamma(\sigma)&=\int_0^\infty s^{\sigma-1}e^{-s}\,ds
 =\int_0^\infty \bigl(t\lambda_j(K)\bigr)^{\sigma-1}e^{-t\lambda_j(K)}\lambda_j(K)\,dt\\
 &=\lambda_j(K)^\sigma\int_0^\infty t^{\sigma-1}e^{-t\lambda_j(K)}\,dt.
\end{align*}
Thus
\[
 \lambda_j(K)^{-\sigma}=\frac{1}{\Gamma(\sigma)}\int_0^\infty t^{\sigma-1}e^{-t\lambda_j(K)}\,dt.
\]
For $\sigma>\frac{n}{2}$, \eqref{eq:heat-trace-upper-bound} and \eqref{eq:heat-trace-exponential-bound} imply $\int_0^\infty t^{\sigma-1}Z_K(t)\,dt<\infty$. Summing over $j$ and applying monotone convergence to the nonnegative integrands, we obtain from \eqref{eq:heat-trace-expansion} that
\begin{equation}\label{eq:zeta-heat-representation}
\begin{aligned}
 \zeta_\sigma(K)&=\sum_{j=1}^\infty\lambda_j(K)^{-\sigma}=\frac{1}{\Gamma(\sigma)}\sum_{j=1}^\infty\int_0^\infty t^{\sigma-1}e^{-t\lambda_j(K)}\,dt\\
 &=\frac{1}{\Gamma(\sigma)}\int_0^\infty t^{\sigma-1}Z_K(t)\,dt<\infty.
\end{aligned}
\end{equation}

Inequality~\eqref{eq:target-zeta-logbm} now follows from Theorem~\ref{thm:transfer-inequalities}~(i), with $H_K=Z_K$ and $\alpha=\sigma$.
\vskip 3pt
\textbf{Case 3.} Prove the inequality \eqref{eq:target-torsion-logbm} for torsional rigidity.

Since $-\Delta u_K=1$ and $-\Delta\phi_j=\lambda_j(K)\phi_j$ hold in the weak sense, it follows that
\[
 \lambda_j(K)\int_Ku_K\phi_j\,dx =\int_K\nabla u_K\cdot\nabla\phi_j\,dx =\int_K\phi_j\,dx.
\]

The preceding equality and Parseval's identity for $u_K,1\in L^2(K)$ imply
\begin{equation}\label{eq:torsion-spectral-expansion}
 \tau(K)=\int_Ku_K\,dx=\sum_{j=1}^\infty\Bigl(\int_Ku_K\phi_j\,dx\Bigr)\Bigl(\int_K\phi_j\,dx\Bigr)=\sum_{j=1}^\infty\frac{\bigl(\int_K\phi_j\,dx\bigr)^2}{\lambda_j(K)}.
\end{equation}
Substituting $\lambda_j(K)^{-1}=\int_0^\infty e^{-t\lambda_j(K)}\,dt$ and using monotone convergence and \eqref{eq:heat-content-expansion}, we conclude that
\begin{equation}\label{eq:torsion-heat-representation}
 \tau(K)=\sum_{j=1}^\infty\int_0^\infty e^{-t\lambda_j(K)}\Bigl(\int_K\phi_j(x)\,dx\Bigr)^2\,dt=\int_0^\infty Q_K(t)\,dt<\infty.
\end{equation}
Inequality~\eqref{eq:target-torsion-logbm} now follows from Theorem~\ref{thm:transfer-inequalities}~(i), with $H_K=Q_K$ and $\alpha=1$.
\vskip 3pt
\textbf{Case 4.} Prove the inequality \eqref{eq:target-chord-logbm} for chord integrals.

Write $a=n+1-q$. Then $0<a<n$ and $E_a(K)<\infty$. By \eqref{ri:eq:chord}, inequality \eqref{eq:target-chord-logbm} is equivalent to the logarithmic Brunn--Minkowski inequality for $E_a$.

The substitution $s=\frac{|x-y|^2}{4t}$ in the gamma integral, together with Fubini's theorem and \eqref{eq:gaussian-integral-definition}, yields
\begin{equation}\label{ri:eq:mellin}
\begin{aligned}
 E_a(K)&=\int_{K^2}|x-y|^{-a}\,dx\,dy\\
 &=\frac1{\Gamma(\frac{a}{2})} \int_{K^2}|x-y|^{-a}\int_0^\infty s^{\frac{a}{2}-1}e^{-s}\,ds\,dx\,dy\\
 &=\frac1{4^{\frac{a}{2}}\Gamma(\frac{a}{2})} \int_{K^2}\int_0^\infty t^{-\frac{a}{2}-1}e^{-\frac{|x-y|^2}{4t}}\,dt\,dx\,dy\\
 &=\frac1{4^{\frac{a}{2}}\Gamma(\frac{a}{2})} \int_0^\infty t^{-\frac{a}{2}-1}F_K(t)\,dt.
\end{aligned}
\end{equation}
In particular, $t^{-\frac{a}{2}-1}F_K(t)$ is integrable on $(0,\infty)$.

By \eqref{eq:gaussian-integral-definition}, $0<F_K(t)\leqslant |K|^2$ for $t>0$. The substitution $x=\sqrt t\,u$, $y=\sqrt t\,v$ gives
\[
 F_K(t)=t^n\int_{(K/\sqrt{t})^2}e^{-\frac{|u-v|^2}{4}}\,du\,dv.
\]
Take $m_j=2$ and $\psi_j(u,v)=\frac{|u-v|^2}{4}$ for $j\geqslant1$. Since $\psi_j$ is even and convex, $F_K$ has representation \eqref{eq:transfer-profile-representation} with $\beta=n$. By \eqref{ri:eq:mellin}, inequality \eqref{eq:target-chord-logbm} now follows from Theorem~\ref{thm:transfer-inequalities}~(i), with $H_K=F_K$ and $\alpha=-\frac{a}{2}$. This completes the proof.
\end{proof}

\vskip 25pt
\section{\bf Equality conditions}\label{sec:rigidity-common}
\vskip 5pt

In this section, we establish the equality conditions in Theorems~\ref{thm:logbm-transfer} and~\ref{thm:geometric-logbm}. Throughout this section, let $K,L\in\Ke$, $0<\theta<1$. Write $M=\comb K L\theta$.

\subsection{Equality conditions in the transfer principle}

\begin{theorem}\label{thm:transfer-equality}
Equality holds in Theorem~\ref{thm:transfer-inequalities}~(i) if and only if, for some $r,s>0$,
\begin{equation}\label{eq:abstract-profile-equality}
 H_L(rt)=sH_K(t),\qquad H_M(r^\theta t)=s^\theta H_K(t),\quad \text{for } t>0.
\end{equation}
\end{theorem}
\begin{proof}
For sufficiency, assume \eqref{eq:abstract-profile-equality} hold. By the definition of $\mathcal J_\alpha$, we have
\begin{align*}
 \mathcal J_\alpha(L)&=\int_0^\infty x^{\alpha-1}H_L(x)\,dx=\int_0^\infty (rt)^{\alpha-1}H_L(rt)\,r\,dt=r^\alpha s\,\mathcal J_\alpha(K),\\
 \mathcal J_\alpha(M)
 &=\int_0^\infty x^{\alpha-1}H_M(x)\,dx
 =\int_0^\infty (r^\theta t)^{\alpha-1}H_M(r^\theta t)\,r^\theta\,dt=(r^\alpha s)^\theta\mathcal J_\alpha(K).
\end{align*}
Hence $\mathcal J_\alpha(M)=\mathcal J_\alpha(K)^{1-\theta}\mathcal J_\alpha(L)^\theta$, as desired.

For necessity, assume equality holds in \eqref{eq:transfer-integral-logbm}.

For $D=K,L,M$, define
\begin{equation}\label{eq:normalized-profile-definition}
 f_D(x)=\frac{e^{\alpha x}H_D(e^x)}{\mathcal J_\alpha(D)},\qquad \text{for }x\in\R.
\end{equation}
By \eqref{eq:abstract-joint-time} with $K=L=D$, $f_D$ is log-concave on $\R$, and hence continuous. Moreover,
\[
 \int_\R f_D(x)\,dx=\frac{1}{\mathcal J_\alpha(D)}\int_\R e^{\alpha x}H_D(e^x)\,dx=\frac{1}{\mathcal J_\alpha(D)}\int_0^\infty t^{\alpha-1}H_D(t)\,dt=1.
\]

Define $T:\R\to\R$ by
\begin{equation}\label{eq:transport-map-definition}
 \int_{-\infty}^{x}f_K(y)\,dy =\int_{-\infty}^{T(x)}f_L(y)\,dy.
\end{equation}
Since $f_K$ and $f_L$ are positive, continuous, and $\int_{\R}f_K(x)\,dx=\int_{\R}f_L(x)\,dx=1$, it follows that $T$ is well defined, $C^1$, and strictly increasing. For $x\in\R$, differentiating \eqref{eq:transport-map-definition} gives
\begin{equation}\label{eq:normalized-derivative-identity}
 f_K(x)=f_L(T(x))T'(x),\qquad T'(x)>0.
\end{equation}

Take $t_0=e^x$ and $t_1=e^{T(x)}$ in \eqref{eq:abstract-joint-time}. Since $\mathcal J_\alpha(M)=\mathcal J_\alpha(K)^{1-\theta}\mathcal J_\alpha(L)^\theta$, it follows from \eqref{eq:normalized-profile-definition} and \eqref{eq:normalized-derivative-identity} that
\begin{equation}\label{eq:normalized-pointwise-bound}
 f_M((1-\theta)x+\theta T(x))\geqslant f_K(x)^{1-\theta}f_L(T(x))^\theta=\frac{f_K(x)}{T'(x)^\theta}.
\end{equation}
By \eqref{eq:normalized-derivative-identity}, the map $x\mapsto(1-\theta)x+\theta T(x)$ is strictly increasing and tends to $\pm\infty$ as $x\to\pm\infty$. Making this change of variables in $\int_\R f_M(x)\,dx=1$ and applying \eqref{eq:normalized-pointwise-bound} and the arithmetic--geometric mean inequality, we obtain
\begin{align*}
 1&=\int_\R f_M((1-\theta)x+\theta T(x)) \bigl(1-\theta+\theta T'(x)\bigr)\,dx\\
 &\geqslant\int_\R f_K(x) \frac{1-\theta+\theta T'(x)}{T'(x)^\theta}\,dx \geqslant\int_\R f_K(x)\,dx=1.
\end{align*}

Therefore, both inequalities are equalities. Since $f_K>0$ and $T'$ is continuous, the equality condition in the arithmetic--geometric mean inequality implies $T'(x)=1$ for $x\in\R$. Thus $T(x)=x+\log r$ for some $r>0$.

Since $T'=1$, both sides of \eqref{eq:normalized-pointwise-bound} integrate to one. By continuity, equality holds for $x\in\R$. Substituting $T(x)=x+\log r$ into \eqref{eq:normalized-derivative-identity} and \eqref{eq:normalized-pointwise-bound}, we obtain
\[
 f_L(x+\log r)=f_K(x),\qquad f_M(x+\theta\log r)=f_K(x).
\]
These, together with \eqref{eq:normalized-profile-definition}, yield \eqref{eq:abstract-profile-equality}, with $s=r^{-\alpha}\frac{\mathcal J_\alpha(L)}{\mathcal J_\alpha(K)}$. This completes the proof.
\end{proof}

\subsection{Rigidity lemmas}

\begin{lemma}\label{lem:measure-support-rigidity}
For $D=K,L,M$, let $\mu_D$ be an even probability measure on $\R^n$ with support $D$. If
\begin{equation}\label{eq:convex-test-inequality}
 \int_{\R^n}\varphi\,d\mu_M \leqslant(1-\theta)\int_{\R^n}\varphi\,d\mu_K +\theta\int_{\R^n}\varphi\,d\mu_L
\end{equation}
for even convex functions $\varphi:\R^n\to\R$, with equality for $\varphi(x)=|x|^2$, then $K=L=M$.
\end{lemma}
\begin{proof}
Let $f\in C_c^\infty(\R^n)$ be even. Since $f$ has bounded second derivatives, we can choose $b>0$ so that $b|x|^2\pm f(x)$ are convex. Apply \eqref{eq:convex-test-inequality} to these two functions. Since equality holds for $\varphi(x)=|x|^2$, it follows that
\[
 \int_{\R^n}f\,d\mu_M=(1-\theta)\int_{\R^n}f\,d\mu_K +\theta\int_{\R^n}f\,d\mu_L.
\]
Thus $\mu_M=(1-\theta)\mu_K+\theta\mu_L$.

Taking supports, we have $K\cup L\subseteq M$. This, together with $M=K^{1-\theta}L^\theta$, yields
\[
 \max\{h_K,h_L\}\leqslant h_M\leqslant h_K^{1-\theta}h_L^\theta \leqslant\max\{h_K,h_L\}.
\]
Equality therefore holds throughout. Thus $h_K=h_L=h_M$ and hence $K=L=M$.
\end{proof}

We use the next lemma to establish the equality conditions for the spectral and torsional functionals in Section~\ref{sec:equality-proof}.

\begin{lemma}\label{lem:ground-state-rigidity}
The first Dirichlet eigenvalue satisfies
\begin{equation}\label{sp:eq:eigenvalue-logbm}
 \lambda_1(M)\leqslant(1-\theta)\lambda_1(K) +\theta\lambda_1(L),
\end{equation}
with equality if and only if $K=L$.
\end{lemma}
\begin{proof}
By \eqref{eq:target-eigen-logbm} and the arithmetic--geometric mean inequality, we have
\[
 \lambda_1(M)\leqslant\lambda_1(K)^{1-\theta}\lambda_1(L)^\theta \leqslant(1-\theta)\lambda_1(K)+\theta\lambda_1(L).
\]
It remains to characterize equality. For $D=K,L,M$, let $k_t^D$ be the Dirichlet heat kernel on $\intt D$, as in Section~\ref{sec:heat-kernels}. Let $\phi_D$ be the positive eigenfunction corresponding to $\lambda_1(D)$, normalized by $\int_D\phi_D^2\,dx=1$, and extend it by zero to $\R^n$.

Fix an even $f\in C_c^\infty(\R^n)$. In the following, we compare the integrals $\int_D f\phi_D^2\,dx$ for $D=K,L,M$. To this end, choose $b>0$ such that $b|x|^2\pm f(x)$ is convex, and define
\begin{equation}\label{sp:eq:U-definition}
 U(x,y)=bx\cdot y+\frac{f(x)+f(y)}2, \quad\text{for }x,y\in\R^n.
\end{equation}
Fix $0<\varepsilon<\frac{1}{2b}$ and define
\[
 g_D(t)=\int_{D^2}k_t^D(x,y)^2 e^{-\frac{\varepsilon U(x,y)}{t}}\,dx\,dy, \quad \text{for }t>0.
\]

\textbf{Step 1.} Prove that
\begin{equation}\label{sp:eq:loop-comparison}
 g_M(t)\geqslant g_K(t)^{1-\theta}g_L(t)^\theta.
\end{equation}

We first approximate $g_D$ using $k_{t,j}^D$ from \eqref{eq:heat-kernel-product}. For $j\geqslant2$, define
\[
 g_{D,j}(t)=\int_{D^2}k_{t,j}^D(x,y)k_{t,j}^D(y,x) e^{-\frac{\varepsilon U(x,y)}{t}}\,dx\,dy, \quad \text{for }t>0.
\]
By Lemma~\ref{lem:heat-kernel-product}, the symmetry of $k_t^D$, and $|\partial D|=0$, we have $k_{t,j}^D(x,y)k_{t,j}^D(y,x)\to k_t^D(x,y)^2$ almost everywhere on $D^2$. Lemma~\ref{lem:heat-kernel-product} also implies
\[
 0\leqslant k_{t,j}^D(x,y)k_{t,j}^D(y,x) \leqslant p_t(x,y)^2.
\]
Since $D$ is compact and $U$ is continuous, $p_t(x,y)^2 e^{-\frac{\varepsilon U(x,y)}{t}}$ is integrable on $D^2$ for fixed $t>0$. By the dominated convergence theorem,
\[
 \lim_{j\to\infty}g_{D,j}(t)=\int_{D^2}k_t^D(x,y)^2 e^{-\frac{\varepsilon U(x,y)}{t}}\,dx\,dy=g_D(t).
\]

Substituting \eqref{eq:heat-kernel-product} into $g_{D,j}$, we have
\[
 g_{D,j}(t)=\Bigl(\frac{4\pi t}{j}\Bigr)^{-nj} \int_{D^{2j}}e^{-\Phi(z)}\,dz,
\]
where $z=(z_0,\ldots,z_{2j-1})$, $z_0=x$, $z_j=y$, $z_{2j}=z_0$, and
\[
 \Phi(z)=\frac j{4t}\sum_{i=1}^{2j}|z_i-z_{i-1}|^2 +\frac{\varepsilon}{t}U(z_0,z_j).
\]

Since $f$ is even, so is $\Phi$. To apply \eqref{eq:weighted-product-comparison}, it remains to prove convexity. Rewrite $\Phi(z)$ as
\begin{align*}
 \Phi(z)&=\Bigl( \frac{j}{4t}\sum_{i=1}^{2j}|z_i-z_{i-1}|^2 -\frac{\varepsilon b}{2t}|z_j-z_0|^2 \Bigr)\\
 &\quad+\frac{\varepsilon}{2t} \bigl(b|z_0|^2+f(z_0)+b|z_j|^2+f(z_j)\bigr).
\end{align*}
The second term is convex by the choice of $b$. For the first term, the Cauchy--Schwarz inequality gives $|z_j-z_0|^2\leqslant j\sum_{i=1}^{2j}|z_i-z_{i-1}|^2$. Since $2\varepsilon b<1$, we have
\[
 \frac{j}{4t}\sum_{i=1}^{2j}|z_i-z_{i-1}|^2 -\frac{\varepsilon b}{2t}|z_j-z_0|^2 \geqslant\frac{j(1-2\varepsilon b)}{4t} \sum_{i=1}^{2j}|z_i-z_{i-1}|^2 \geqslant0.
\]
The first term is therefore a nonnegative quadratic form, so $\Phi$ is convex. Applying \eqref{eq:weighted-product-comparison} with $m=2j$ and $\psi=\Phi$, we obtain
\[
 g_{M,j}(t)\geqslant g_{K,j}(t)^{1-\theta}g_{L,j}(t)^\theta.
\]
Letting $j\to\infty$ yields \eqref{sp:eq:loop-comparison}.
\vskip 3pt
\textbf{Step 2.} Prove that
\begin{equation}\label{sp:eq:loop-log-asymptotics}
 \log g_D(t)=-2t\lambda_1(D) -\frac{\varepsilon}{t}\int_D f\phi_D^2\,dx+O(t^{-2}) \quad\text{as }t\to\infty.
\end{equation}

Let $q_t^D=e^{t\lambda_1(D)}k_t^D$ and $\delta_D=\lambda_2(D)-\lambda_1(D)$. Then $\delta_D>0$. Since $k_1^D\leqslant p_1\leqslant(4\pi)^{-\frac n2}$ and $D$ is bound, we have $q_1^D\in L^2(D^2)$. For $t\geqslant1$, it follows from \eqref{eq:heat-kernel-spectral} and the orthonormality of the eigenfunctions that
\begin{align*}
 &\|q_t^D-\phi_D\otimes\phi_D\|_{L^2(D^2)}^2 =\sum_{k=2}^\infty e^{-2t(\lambda_k(D)-\lambda_1(D))}\\
 \leqslant{}& e^{-2(t-1)\delta_D}\|q_1^D\|_{L^2(D^2)}^2 =O(e^{-2\delta_Dt}),
\end{align*}
where $(\phi_D\otimes\phi_D)(x,y)=\phi_D(x)\phi_D(y)$.

Since $g_D(t)$ involves the square of the kernel, we next estimate the difference of the squared kernels in $L^1(D^2)$. Observe that $\|\phi_D\otimes\phi_D\|_{L^2(D^2)}=\int_D\phi_D^2\,dx=1$, so the triangle inequality implies $\|q_t^D\|_{L^2(D^2)}=O(1)$. By the Cauchy--Schwarz inequality,
\begin{equation}\label{sp:eq:kernel-square-L1}
\begin{aligned}
 &\|(q_t^D)^2-(\phi_D\otimes\phi_D)^2\|_{L^1(D^2)}\\
 \leqslant{}& \|q_t^D-\phi_D\otimes\phi_D\|_{L^2(D^2)} \bigl(\|q_t^D\|_{L^2(D^2)}+1\bigr) =O(e^{-\delta_Dt}).
\end{aligned}
\end{equation}

Since $-D=D$, the function $x\mapsto\phi_D(-x)$ is also a positive normalized eigenfunction corresponding to $\lambda_1(D)$. This eigenvalue is simple, so $\phi_D(-x)=\phi_D(x)$. Thus $\phi_D^2\,dx$ is an even probability measure and $\int_Dx\phi_D(x)^2\,dx=0$. By \eqref{sp:eq:U-definition} and $\int_D\phi_D^2\,dx=1$, we have
\[
 \int_{D^2}U(x,y)\phi_D(x)^2\phi_D(y)^2\,dx\,dy =b\Bigl|\int_Dx\phi_D(x)^2\,dx\Bigr|^2 +\int_D f\phi_D^2\,dx=\int_D f\phi_D^2\,dx.
\]
Since $U$ is bounded on $D^2$, Taylor's formula gives $e^{-\frac{\varepsilon U}{t}}=1-\frac{\varepsilon U}{t}+O(t^{-2})$ uniformly on $D^2$ as $t\to\infty$. Using \eqref{sp:eq:kernel-square-L1} and $e^{-\delta_Dt}=O(t^{-2})$, we obtain
\begin{align*}
 e^{2t\lambda_1(D)}g_D(t) &=\int_{D^2}(q_t^D(x,y))^2 e^{-\frac{\varepsilon U(x,y)}{t}}\,dx\,dy\\
 &=\int_{D^2}\Bigl(1-\frac{\varepsilon}{t}U(x,y)\Bigr) \phi_D(x)^2\phi_D(y)^2\,dx\,dy+O(t^{-2})\\
 &=1-\frac{\varepsilon}{t}\int_D f\phi_D^2\,dx+O(t^{-2}),
\end{align*}
which yields \eqref{sp:eq:loop-log-asymptotics}.
\vskip 3pt

\textbf{Step 3.} Characterize equality.

Assume equality in \eqref{sp:eq:eigenvalue-logbm}. The terms of order $t$ in \eqref{sp:eq:loop-log-asymptotics} then cancel, so \eqref{sp:eq:loop-comparison} gives
\begin{align*}
 0&\leqslant \log g_M(t)-(1-\theta)\log g_K(t)-\theta\log g_L(t)\\
 &=\frac{\varepsilon}{t} \Bigl((1-\theta)\int_K f\phi_K^2\,dx +\theta\int_L f\phi_L^2\,dx-\int_M f\phi_M^2\,dx\Bigr) +O(t^{-2}).
\end{align*}
Multiplying by $\frac{t}{\varepsilon}$ and letting $t\to\infty$ yields
\[
 \int_M f\phi_M^2\,dx\leqslant(1-\theta)\int_K f\phi_K^2\,dx +\theta\int_L f\phi_L^2\,dx.
\]
Applying this inequality to $-f$, we obtain equality for even $f\in C_c^\infty(\R^n)$. Since the measures $\phi_D^2\,dx$ are even for $D=K,L,M$, it follows that
\begin{equation}\label{sp:eq:ground-state-measure-identity}
 \phi_M^2\,dx=(1-\theta)\phi_K^2\,dx+\theta\phi_L^2\,dx.
\end{equation}

The function $\phi_D$ is positive in $\intt D$ and vanishes outside $D$, so the measure $\phi_D^2\,dx$ has support $D$. By \eqref{sp:eq:ground-state-measure-identity}, equality holds in the integral inequality of Lemma~\ref{lem:measure-support-rigidity} for even convex functions, so $K=L=M$.

Conversely, equality holds in \eqref{sp:eq:eigenvalue-logbm} when $K=L$. This completes the proof.
\end{proof}

Recall that $F_K$ is defined by \eqref{eq:gaussian-integral-definition}.

\begin{lemma}\label{ri:lem:rigidity}
If $F_M(t)=F_K(t)^{1-\theta}F_L(t)^\theta$ for $t>0$, then $K=L=M$.
\end{lemma}
\begin{proof}
For $D=K,L,M$, let $\mu_D$ be the probability measure on $\R^n$ given by $d\mu_D(x)=\frac{\ind_D(x)}{|D|}\,dx$. Its support is $D$, and it is even since $D=-D$. It remains to verify the conditions in Lemma~\ref{lem:measure-support-rigidity}.

Since $D$ is bounded, we have
\[
 \exp\!\Bigl(-\frac{|x-y|^2}{4t}\Bigr)=1-\frac{|x-y|^2}{4t}+O(t^{-2})
\]
uniformly on $D^2$ as $t\to\infty$. Since $\int_D x\,dx=0$ by symmetry,
\[
 \int_D\int_D |x-y|^2\,dy\,dx=2|D|\int_D|x|^2\,dx-2\Bigl|\int_D x\,dx\Bigr|^2=2|D|\int_D|x|^2\,dx.
\]
It follows from \eqref{eq:gaussian-integral-definition} that
\begin{align*}
 \log F_D(t) &=2\log|D|-\frac{1}{4t|D|^2}\int_{D^2}|x-y|^2\,dx\,dy+O(t^{-2})\\
 &=2\log|D|-\frac{1}{2t}\int_D|x|^2\,d\mu_D+O(t^{-2}).
\end{align*}
Taking logarithms in $F_M(t)=F_K(t)^{1-\theta}F_L(t)^\theta$ and substituting this expansion, we obtain
\[
 2\log\frac{|M|}{|K|^{1-\theta}|L|^\theta}=\frac{1}{2t}\Bigl(\int_M|x|^2\,d\mu_M-(1-\theta)\int_K|x|^2\,d\mu_K-\theta\int_L|x|^2\,d\mu_L\Bigr)+O(t^{-2}).
\]
Letting $t\to\infty$ yields $|M|=|K|^{1-\theta}|L|^\theta$. Multiplying by $2t$ and letting $t\to\infty$ again yields
\begin{equation}\label{ri:eq:second-moment-identity}
 \int_M|x|^2\,d\mu_M =(1-\theta)\int_K|x|^2\,d\mu_K +\theta\int_L|x|^2\,d\mu_L.
\end{equation}

Let $\varphi:\R^n\to\R$ be even and convex. Applying Lemma~\ref{lem:saroglou-density} with $\psi=r\varphi$ gives
\[
 \int_M e^{-r\varphi}\,d\mu_M\geqslant \Bigl(\int_K e^{-r\varphi}\,d\mu_K\Bigr)^{1-\theta} \Bigl(\int_L e^{-r\varphi}\,d\mu_L\Bigr)^\theta,\quad\text{for } r\geqslant0.
\]
Both sides equal $1$ at $r=0$. Since $\varphi$ is bounded on $K\cup L\cup M$, we may compare their right derivatives there to obtain
\[
 \int_M \varphi\,d\mu_M\leqslant (1-\theta)\int_K \varphi\,d\mu_K +\theta\int_L \varphi\,d\mu_L.
\]
Applying this and \eqref{ri:eq:second-moment-identity} in Lemma~\ref{lem:measure-support-rigidity} yields $K=L=M$.
\end{proof}

\subsection{Equality conditions for geometric functionals}\label{sec:equality-proof}

\begin{theorem}
In each inequality of Theorem~\ref{thm:functional-inequalities}, equality holds if and only if $K$ and $L$ are dilates.
\end{theorem}
\begin{proof}
If $L=tK$ for some $t>0$, then $M=t^\theta K$. By homogeneity, all four inequalities in Theorem~\ref{thm:functional-inequalities} hold with equality. It remains to prove the converse.

\textbf{Case 1.} Assume that equality holds in \eqref{eq:target-eigen-logbm} for the first Dirichlet eigenvalue $\lambda_1$. Write $\widetilde{K}=\sqrt{\lambda_1(K)}K$, $\widetilde{L}=\sqrt{\lambda_1(L)}L$. By homogeneity, we have $\lambda_1(\widetilde{K})=\lambda_1(\widetilde{L})=1$. This, together with the assumption, yields
\[
 \lambda_1(\comb{\widetilde{K}}{\widetilde{L}}\theta)=\lambda_1(\widetilde{K})^{1-\theta}\lambda_1(\widetilde{L})^\theta=1.
\]
Lemma~\ref{lem:ground-state-rigidity} therefore implies $\widetilde{K}=\widetilde{L}$. Thus $K$ and $L$ are dilates.

\textbf{Case 2.} Assume that equality holds in \eqref{eq:target-zeta-logbm} for the spectral zeta function $\zeta_\sigma$ with $\sigma>\frac{n}{2}$, or in \eqref{eq:target-torsion-logbm} for torsional rigidity $\tau$. By \eqref{eq:zeta-heat-representation} and \eqref{eq:torsion-heat-representation}, equality holds in \eqref{eq:transfer-integral-logbm} with $H_K=Z_K$ and $\alpha=\sigma$ in the spectral zeta case, and with $H_K=Q_K$ and $\alpha=1$ in the torsional case. Theorem~\ref{thm:transfer-equality} now yields \eqref{eq:abstract-profile-equality} for some $r,s>0$. Applying \eqref{eq:heat-profile-decay} to \eqref{eq:abstract-profile-equality}, we obtain $\lambda_1(L)=r^{-1}\lambda_1(K)$ and
\[
 \lambda_1(M)=r^{-\theta}\lambda_1(K)=\lambda_1(K)^{1-\theta}\lambda_1(L)^\theta.
\]
By Case 1, $K$ and $L$ are dilates.

\textbf{Case 3.} Assume that equality holds in \eqref{eq:target-chord-logbm} for the chord integral $I_q$. By \eqref{ri:eq:chord}, this is equivalent to equality for $E_a$, where $a=n+1-q\in(0,n)$. The representation \eqref{ri:eq:mellin} then implies equality in \eqref{eq:transfer-integral-logbm} with $H_K=F_K$ and $\alpha=-\frac{a}{2}$. Theorem~\ref{thm:transfer-equality} now yields \eqref{eq:abstract-profile-equality} for some $r,s>0$.

Write $L_0=r^{-\frac{1}{2}}L$ and $M_0=r^{-\frac{\theta}{2}}M=\comb K{L_0}\theta$.
By a change of variables in \eqref{eq:gaussian-integral-definition}, $F_{cD}(t)=c^{2n}F_D\bigl(\frac{t}{c^2}\bigr)$ for $c>0$. Thus the identities in \eqref{eq:abstract-profile-equality} take the form
\[
 F_{L_0}(t)=sr^{-n}\,F_K(t),\qquad F_{M_0}(t)=(sr^{-n})^\theta F_K(t).
\]
Consequently, $F_{M_0}(t)=F_K(t)^{1-\theta}F_{L_0}(t)^\theta$
for $t>0$. Lemma~\ref{ri:lem:rigidity} now implies $L_0=K$, and hence $L=\sqrt{r}\,K$.
\end{proof}

Together with Theorem~\ref{thm:functional-inequalities}, this completes the proof of Theorem~\ref{thm:geometric-logbm}.

\vskip 25pt
\section{\bf A variational principle for logarithmic Minkowski problems}\label{sec:abstract-inverse}
\vskip 5pt

In this section, we establish a variational principle for the logarithmic Minkowski problems studied in Section~\ref{sec:inverse-applications}.

Let $\mathcal F:\Ke\to(0,\infty)$ be monotone with respect to inclusion and homogeneous of degree one. Throughout this section, assume that the following variation formula holds.

\begin{itemize}
\item[\textbf{(V)}]\phantomsection\label{cond:variation} For $K\in\Ke$, there exists a finite even Borel measure $\eta_K$ on $\Sph$ such that
\begin{equation*}
 \frac d{dt}\log\mathcal F([h_Ke^{tg}])\Bigr|_{t=0} =\int_{\Sph}g\,d\eta_K,\qquad\text{for }g\in C_e(\Sph).
\end{equation*}
\end{itemize}

By \hyperref[cond:variation]{(V)} and homogeneity, $\eta_K$ is a probability measure and $\eta_{rK}=\eta_K$ for $r>0$. To give a unified treatment of the four cases in Theorem~\ref{thm:geometric-inverse}, we study the following problem.
\vskip 3pt
\noindent\textbf{General logarithmic Minkowski problem.} \emph{Given an even Borel probability measure $\mu$ on $\Sph$, what are necessary and sufficient conditions for the existence of $K\in\Ke$ satisfying $\eta_K=\mu$? If such a body exists, is it unique up to dilation?}

If $F:\Ke\to(0,\infty)$ is homogeneous of degree $s\ne0$, we can reduce its logarithmic Minkowski problem to the one above by working with $F^{\frac{1}{s}}$. See Section~\ref{sec:geometric-inverse-proof} for details.

\subsection{Uniqueness}

Uniqueness follows from the following logarithmic Brunn--Minkowski inequality.

\begin{itemize}
\item[\textbf{(L)}]\phantomsection\label{cond:logbm} For $K,L\in\Ke$ and $0<\theta<1$,
\begin{equation}\label{eq:abstract-strict-logbm}
 \mathcal F(\comb K L\theta) \geqslant\mathcal F(K)^{1-\theta}\mathcal F(L)^\theta,
\end{equation}
with equality if and only if $L=rK$ for some $r>0$.
\end{itemize}

For volume, B\"or\"oczky--Lutwak--Yang--Zhang \cite[Lemma~3.2]{BLYZ12} proved the equivalence of the logarithmic Brunn--Minkowski and logarithmic Minkowski inequalities.

\begin{lemma}\label{prop:abstract-measure}
Condition \hyperref[cond:logbm]{(L)} is equivalent to the assertion that for all $K,L\in\Ke$,
\begin{equation}\label{eq:abstract-logmink}
 \int_{\Sph}\log\frac{h_L}{h_K}\,d\eta_K \geqslant\log\frac{\mathcal F(L)}{\mathcal F(K)},
\end{equation}
with equality if and only if $K$ and $L$ are dilates.
\end{lemma}
\begin{proof}
Assume \hyperref[cond:logbm]{(L)} holds. Let $K,L\in\Ke$. Write $K_t=\comb K L t$ for $0\leqslant t\leqslant1$. Then $K_t=[h_K e^{t\log\frac{h_L}{h_K}}]$. By \hyperref[cond:variation]{(V)} and \eqref{eq:abstract-strict-logbm}, we have
\[
 \int_{\Sph}\log\frac{h_L}{h_K}\,d\eta_K=\lim_{t\to 0^+}\frac{\log\mathcal F(K_t)-\log\mathcal F(K)}{t}\geqslant\log\frac{\mathcal F(L)}{\mathcal F(K)}.
\]
Thus \eqref{eq:abstract-logmink} holds.

Assume equality holds in \eqref{eq:abstract-logmink}. Let $0<t<1$. Since $K_t=[h_K^{1-t}h_L^t]$, we have $\bigl(\frac{h_L}{h_K}\bigr)^t\geqslant\frac{h_{K_t}}{h_K}$. By this, \eqref{eq:abstract-logmink} for $K,K_t$, and \eqref{eq:abstract-strict-logbm} for $K,L$, it follows that
\[
 t\int_{\Sph}\log\frac{h_L}{h_K}\,d\eta_K \geqslant\int_{\Sph}\log\frac{h_{K_t}}{h_K}\,d\eta_K \geqslant\log\frac{\mathcal F(K_t)}{\mathcal F(K)} \geqslant t\log\frac{\mathcal F(L)}{\mathcal F(K)}.
\]
Since the first and last terms agree by assumption, \eqref{eq:abstract-strict-logbm} holds with equality, and hence $K$ and $L$ are dilates. Conversely, if $L=rK$ for some $r>0$, both sides of \eqref{eq:abstract-logmink} equal $\log r$ since $\eta_K(\Sph)=1$ and $\mathcal F(rK)=r\mathcal F(K)$.

To prove the reverse implication, assume \eqref{eq:abstract-logmink} holds for all $K,L\in\Ke$, with equality only for dilates. Let $K,L\in\Ke$ and $0<\theta<1$. Write $M=\comb K L\theta$ and $g=\log\frac{h_K^{1-\theta}h_L^\theta}{h_M}$. Then for $0\leqslant s\leqslant1$, $h_M\leqslant h_M e^{sg}\leqslant h_K^{1-\theta}h_L^\theta$. Since $[h_M]=M=[h_K^{1-\theta}h_L^\theta]$, it follows that $[h_M e^{sg}]=M$. Hence $\log\mathcal F([h_M e^{sg}])$ is constant for $0\leqslant s\leqslant1$. By \hyperref[cond:variation]{(V)}, we have
\[
 \int_{\Sph}g\,d\eta_M=\frac{d}{ds}\log\mathcal F([h_M e^{sg}])\Bigr|_{s=0^+}=0.
\]

From $g=\log\frac{h_K^{1-\theta}h_L^\theta}{h_M}$ and \eqref{eq:abstract-logmink}, it follows that
\begin{align*}
 0=\int_{\Sph}g\,d\eta_M
 &=(1-\theta)\int_{\Sph}\log\frac{h_K}{h_M}\,d\eta_M+\theta\int_{\Sph}\log\frac{h_L}{h_M}\,d\eta_M\\
 &\geqslant(1-\theta)\log\frac{\mathcal F(K)}{\mathcal F(M)}+\theta\log\frac{\mathcal F(L)}{\mathcal F(M)},
\end{align*}
which implies \eqref{eq:abstract-strict-logbm}.

Assume equality holds in \eqref{eq:abstract-strict-logbm}. Then equality must hold in both applications of \eqref{eq:abstract-logmink}. Thus $K$ and $L$ are both dilates of $M$. Conversely, if $L=rK$, then $M=[r^\theta h_K]=r^\theta K$. By homogeneity, equality holds in \eqref{eq:abstract-strict-logbm}.
\end{proof}

\begin{theorem}\label{thm:abstract-uniqueness}
Suppose that $\mathcal F$ satisfies \hyperref[cond:logbm]{(L)}. Then for $K,L\in\Ke$, $\eta_K=\eta_L$ if and only if $K$ and $L$ are dilates.
\end{theorem}
\begin{proof}
Assume $\eta_K=\eta_L$. By \eqref{eq:abstract-logmink}, we have
\[
 \log\frac{\mathcal F(L)}{\mathcal F(K)} \leqslant\int_{\Sph}\log\frac{h_L}{h_K}\,d\eta_K =-\int_{\Sph}\log\frac{h_K}{h_L}\,d\eta_L \leqslant\log\frac{\mathcal F(L)}{\mathcal F(K)}.
\]
Thus, both inequalities are equalities. By Lemma~\ref{prop:abstract-measure}, $K$ and $L$ are dilates. Conversely, if $L=rK$ for some $r>0$, then $\eta_L=\eta_K$, as noted after \hyperref[cond:variation]{(V)}.
\end{proof}

\subsection{Existence}

Fix $b_1>0$ and $b_2,\ldots,b_n\geqslant0$ with $\sum_{i=1}^n b_i=1$, and write $c_k=\sum_{i=1}^k b_i$ for $1\leqslant k<n$. We use the following conditions involving cylinders and John ellipsoids.

\begin{itemize}
\item[\textbf{(C)}]\phantomsection\label{cond:cylinder} For every $k$-dimensional subspace $\xi$ of $\R^n$, $1\leqslant k<n$,
\begin{equation*}
 \lim_{r\to\infty} \frac{\log\mathcal F(B_\xi\times rB_{\xi^\perp})}{\log r}=1-c_k.
\end{equation*}
\item[\textbf{(J)}]\phantomsection\label{cond:john} For $\varepsilon>0$, there exists $C_\varepsilon>0$ such that, for $K\in\Ke$,
\begin{equation*}
 \mathcal F(K)\leqslant C_\varepsilon\Bigl(\frac{a_n}{a_1}\Bigr)^\varepsilon\prod_{i=1}^n a_i^{b_i},
\end{equation*}
where $0<a_1\leqslant\cdots\leqslant a_n$ are the semiaxis lengths of the John ellipsoid $\mathrm{J}K$, the unique ellipsoid of maximal volume contained in $K$.
\end{itemize}

Let $\mathcal P_c$ be the class of even probability measures $\mu$ on $\Sph$ satisfying $\mu(\Sph\cap \xi)<c_k$ for each $k=1,\ldots,n-1$ and every $k$-dimensional subspace $\xi$ of $\R^n$.

\begin{theorem}\label{thm:abstract-inverse-principle}
\leavevmode
\begin{enumerate}
\item[(i)] If $\mathcal{F}$ satisfies \hyperref[cond:logbm]{(L)} and \hyperref[cond:cylinder]{(C)}, then $\eta_K\in\mathcal P_c$ for $K\in\Ke$.
\item[(ii)] If $\mathcal{F}$ satisfies \hyperref[cond:john]{(J)} and $\mu\in\mathcal P_c$, then there exists $K\in\Ke$ such that $\eta_K=\mu$.
\end{enumerate}
\end{theorem}

\begin{proof}
\textbf{(i)} Let $K\in\Ke$ and $\xi$ be a $k$-dimensional subspace with $0<k<n$. Choose $r_0>0$ so that $r_0(B_\xi\times B_{\xi^\perp})\subseteq K$. For $t\geqslant0$, let $L_t=r_0(B_\xi\times e^tB_{\xi^\perp})$. Then
\begin{align*}
 L_t\subseteq e^t r_0(B_\xi\times B_{\xi^\perp})\subseteq e^tK\quad\text{and}\quad \mathrm{P}_{\xi}L_t=r_0B_\xi\subseteq\mathrm{P}_{\xi}K.
\end{align*}

Thus, for $t>1$, $h_K^{1-\frac{1}{t}}h_{L_t}^{\frac{1}{t}}\leqslant h_K^{1-\frac{1}{t}}(e^th_K)^{\frac{1}{t}}=eh_K$ on $\Sph$ and $h_K^{1-\frac{1}{t}}h_{L_t}^{\frac{1}{t}}\leqslant h_K$ on $\Sph\cap\xi$. From the definition of the Wulff shape, it follows that
\[\comb K{L_t}{\frac{1}{t}}=[h_K^{1-\frac{1}{t}}h_{L_t}^{\frac{1}{t}}]\subseteq eK\cap(\mathrm{P}_{\xi}K+\xi^\perp)\eqqcolon D.\]
By \hyperref[cond:logbm]{(L)} and the monotonicity and homogeneity of $\mathcal F$, we have
\[
 \log\frac{\mathcal F(D)}{\mathcal F(K)}\geqslant\frac1t\log\frac{\mathcal F(L_t)}{\mathcal F(K)}=\frac1t\log\frac{r_0\,\mathcal F(B_\xi\times e^tB_{\xi^\perp})}{\mathcal F(K)}.
\]
Letting $t\to\infty$ and using \hyperref[cond:cylinder]{(C)} with $r=e^t$ yields
\begin{equation}\label{eq:cylinder-limit-lower-bound}
 \log\frac{\mathcal F(D)}{\mathcal F(K)}\geqslant 1-c_k.
\end{equation}

Since $K\subseteq D\subseteq\mathrm{P}_{\xi}K+\xi^\perp$, we have $\mathrm{P}_{\xi}D=\mathrm{P}_{\xi}K$. If $D=sK$ for some $s>0$, this identity forces $s=1$, while $D\cap\xi^\perp=e(K\cap\xi^\perp)$ forces $s=e$. Thus $K$ and $D$ are not dilates. By $\eta_K(\Sph)=1$, $h_D\leqslant eh_K$ on $\Sph$, $h_D=h_K$ on $\Sph\cap \xi$, the strict inequality in Lemma~\ref{prop:abstract-measure}, and \eqref{eq:cylinder-limit-lower-bound}, it follows that
\begin{align*}
 1-\eta_K(\Sph\cap\xi) &=\eta_K(\Sph\setminus\xi)\geqslant\int_{\Sph\setminus\xi}\log\frac{h_D}{h_K}\,d\eta_K\\
 &=\int_{\Sph}\log\frac{h_D}{h_K}\,d\eta_K>\log\frac{\mathcal F(D)}{\mathcal F(K)}\geqslant 1-c_k.
\end{align*}
Thus $\eta_K(\Sph\cap \xi)<c_k$, as desired.
\vskip 3pt
\textbf{(ii)} Let $\mu\in\mathcal P_c$. Define
\[
 \Phi_\mu(K)=\int_{\Sph}\log h_K\,d\mu-\log\mathcal F(K),\quad\text{for }K\in\Ke.
\]
We divide the proof of (ii) into three steps.
\vskip 3pt
\textbf{Step 1.} Estimate $\Phi_\mu$.

Write $d(u,\xi)$ for the Euclidean distance from $u$ to $\xi$. We claim that there exist $\delta>0$ and $\varepsilon\in(0,1)$ such that
\begin{equation}\label{eq:tube}
 \mu(\{u\in\Sph:d(u,\xi)\leqslant\varepsilon\})\leqslant c_k-\delta
\end{equation}
for every $k$-dimensional subspace $\xi$ of $\R^n$, $1\leqslant k<n$. Otherwise, by compactness of the Grassmannians, there are $k$-dimensional subspaces $\xi_j\to \xi$, $1\leqslant k<n$, with
\[
 \mu(\{u\in\Sph:d(u,\xi_j)\leqslant\frac1j\})>c_k-\frac1j.
\]
For fixed $r>0$, these shrinking tubes eventually lie in $\{u\in\Sph:d(u,\xi)\leqslant r\}$. Letting $j\to\infty$, we have
\[
 \mu(\{u\in\Sph:d(u,\xi)\leqslant r\})\geqslant c_k.
\]
As $r\to0^+$, these closed neighborhoods decrease to $\Sph\cap \xi$. By continuity from above, $\mu(\Sph\cap \xi)\geqslant c_k$, contrary to $\mu\in\mathcal P_c$.

For $K\in\Ke$, choose orthonormal principal directions $e_i$ of $\mathrm{J}K$ corresponding to $a_i$, and let $\xi_k=\operatorname{span}\{e_1,\ldots,e_k\}$ for $1\leqslant k<n$. Since $\mathrm{J}K\subseteq K$, we have $h_K\geqslant a_1$ and $h_K(u)\geqslant a_{k+1}d(u,\xi_k)$. For $a_k\leqslant r\leqslant a_{k+1}$, \eqref{eq:tube} therefore implies
\[
 \mu(\{u:h_K(u)>\varepsilon r\})\geqslant1-c_k+\delta.
\]
By Fubini's theorem, the preceding estimate, $c_k=\sum_{i=1}^k b_i$ and $\sum_{i=1}^n b_i=1$, we have
\begin{align*}
 \int_{\Sph}\log h_K\,d\mu &=\log(\varepsilon a_1)+\int_{a_1}^{\infty}\mu(\{u:h_K(u)>\varepsilon r\})\,\frac{dr}{r}\\
 &\geqslant\log(\varepsilon a_1)+\sum_{k<n}(1-c_k+\delta)\log\frac{a_{k+1}}{a_k}\\
 &=\sum_i b_i\log a_i+\delta\log\frac{a_n}{a_1}+\log\varepsilon.
\end{align*}
Taking the exponent in \hyperref[cond:john]{(J)} to be $\frac{\delta}{2}$ and applying the preceding estimate, we obtain
\begin{equation*}
 \Phi_\mu(K)\geqslant\frac{\delta}{2}\log\frac{a_n}{a_1}-C,
\end{equation*}
where $C>0$ is independent of $K$. In particular, $\Phi_\mu$ is bounded below.
\vskip 3pt
\textbf{Step 2.} Find a minimizer.

Let $\{K_j\}$ be a minimizing sequence for $\Phi_\mu$. Since $\mu$ is a probability measure and $\mathcal F$ is homogeneous of degree one, $\Phi_\mu(rK)=\Phi_\mu(K)$ for $r>0$. We may therefore assume $a_1(K_j)=1$. Since $\Phi_\mu(K_j)$ is bounded above, the estimate in Step~1 gives $\sup_j a_n(K_j)<\infty$. John's theorem (see, e.g., \cite[Theorem~10.12.2]{Schneider}) then gives
\[
 B^n\subseteq K_j\subseteq RB^n
\]
for a fixed $R$. By the Blaschke selection theorem, $\{K_j\}$ has a convergent subsequence, denoted again by $\{K_j\}$, with limit $K\in\Ke$. From uniform convergence of $\log h_{K_j}$ and Lemma~\ref{lem:wulff-comparison}, it follows that
\[
 \Phi_\mu(K)=\lim_{j\to\infty}\Phi_\mu(K_j)=\inf_{L\in\Ke}\Phi_\mu(L).
\]

\textbf{Step 3.} Prove that $\eta_K=\mu$.

For $g\in C_e(\Sph)$, let $K_t=[h_Ke^{tg}]$. Consider the function
\[
 \varphi(t)=\int_{\Sph}(\log h_K+tg)\,d\mu-\log\mathcal F(K_t).
\]
Since $K$ minimizes $\Phi_\mu$ and $h_{K_t}\leqslant h_Ke^{tg}$, we have
\[
 \varphi(0)=\Phi_\mu(K)\leqslant\Phi_\mu(K_t)\leqslant\varphi(t).
\]
Thus $\varphi$ has a minimum at $t=0$. By \hyperref[cond:variation]{(V)},
\[
 0=\varphi'(0)=\int_{\Sph}g\,d\mu-\int_{\Sph}g\,d\eta_K.
\]
Since $g$ is arbitrary and both measures are even, $\eta_K=\mu$.
\end{proof}

\vskip 25pt

\section{\bf Logarithmic Minkowski problems for geometric functionals}\label{sec:inverse-applications}
\vskip 5pt

In this section, we prove Theorem~\ref{thm:geometric-inverse} by applying Theorems~\ref{thm:abstract-uniqueness} and~\ref{thm:abstract-inverse-principle}. Condition \hyperref[cond:logbm]{(L)} follows from Theorem~\ref{thm:geometric-logbm}, so it remains to verify \hyperref[cond:variation]{(V)}, \hyperref[cond:john]{(J)}, and \hyperref[cond:cylinder]{(C)}.

\subsection{The cone spectral zeta measure}\label{sec:spectral-measures}

Let $K\in\Ko$ and $\{\phi_j\}$ be the Dirichlet eigenbasis as in Section~\ref{sec:heat-kernels}. For $\sigma>\frac{n}{2}$, define the cone spectral zeta measure of $K$ for Borel $\omega\subseteq\Sph$ by
\begin{equation}\label{sp:eq:zetameasure}
 \zeta_\sigma(K,\omega)=\sigma\int_{\nu_K^{-1}(\omega)}h_K(\nu_K(x))\sum_{j=1}^\infty\lambda_j(K)^{-\sigma-1}|\partial_\nu\phi_j(x)|^2\,d\mathcal{H}^{n-1}(x).
\end{equation}
By Grisvard's regularity theorem \cite[Theorem~3.2.1.2]{Grisvard}, $\phi_j\in H^2(\intt K)$, so the normal component $\partial_\nu\phi_j$ of the boundary trace of $\nabla\phi_j$ belongs to $L^2(\partial K)$. Note that $\zeta_\sigma(K,\cdot)$ is absolutely continuous with respect to $S_K$. The sum is independent of the choice of orthonormal eigenbasis. If $K=-K$, reflection preserves each eigenspace, so the measure is even.

Let $u$ be a Dirichlet eigenfunction on $\intt K$ with eigenvalue $\lambda$ and $X$ a $C^1$ vector field on an open neighborhood of $K$. Since $u$ has zero trace, the boundary trace of $\nabla u$ is $(\partial_\nu u)\nu_K$. Substituting $\Delta u=-\lambda u$ in the Rellich formula \cite[p.~223, equation~(1)]{Ancona} and integrating $-\lambda X\cdot\nabla(u^2)$ by parts, we obtain
\begin{equation}\label{sp:eq:zetarellich}
 \int_{\partial K}(X\cdot\nu_K)|\partial_\nu u|^2\,d\mathcal{H}^{n-1} =\int_K\bigl(2\langle DX\nabla u,\nabla u\rangle -(\operatorname{div}X)|\nabla u|^2+\lambda(\operatorname{div}X)u^2\bigr)\,dx,
\end{equation}
where $DX$ is the Jacobian matrix of $X$.

Taking $X(x)=x$ and $u=\phi_j$ in \eqref{sp:eq:zetarellich} gives $\int_{\partial K}h_K(\nu_K)|\partial_\nu\phi_j|^2\,d\mathcal{H}^{n-1}=2\lambda_j(K)$. By \eqref{sp:eq:zetameasure} and monotone convergence,
\[
 \zeta_\sigma(K,\Sph)=2\sigma\zeta_\sigma(K)<\infty.
\]
Since $\zeta_\sigma(K,\Sph)<\infty$ and $\min_{\Sph}h_K>0$, the series in \eqref{sp:eq:zetameasure} converges in $L^1(\partial K)$.

\begin{lemma}\label{sp:prop:zetaweak}
Let $\sigma>\frac{n}{2}$. If $K_i\to K$ in $\Ko$ as $i\to\infty$, then $\zeta_\sigma(K_i,\cdot)\to\zeta_\sigma(K,\cdot)$ weakly.
\end{lemma}

\begin{proof}
Write $\Omega=\intt K$. Let $R>0$ be such that $K_i,K\subseteq RB^n$ for all $i$.
\vskip 3pt
\textbf{Step 1.} Since $K_i\to K$ in $\Ko$, we can choose $\varepsilon_i\to0$ so that $(1-\varepsilon_i)K\subseteq K_i\subseteq(1+\varepsilon_i)K$ for large $i$. By domain monotonicity and homogeneity \eqref{eq:eigenvalue-monotonicity-scaling}, we have
\[
 \frac{\lambda_j(K)}{(1+\varepsilon_i)^2}\leqslant\lambda_j(K_i)\leqslant\frac{\lambda_j(K)}{(1-\varepsilon_i)^2}.
\]
Thus, for fixed $j$, $\lambda_j(K_i)\to\lambda_j(K)$ as $i\to\infty$.

Extend an orthonormal Dirichlet eigenbasis $\{\phi_{j,i}\}_{j\geqslant1}$ of $K_i$ by zero to $RB^n$. For fixed $j$, these extensions are bounded in $H^1$, since $\|\phi_{j,i}\|_{H^1}^2=1+\lambda_j(K_i)\to1+\lambda_j(K)$. By the Rellich--Kondrachov theorem \cite[Theorem~7.22]{GilbargTrudinger} and a diagonal argument, from any subsequence of $(K_i)$ we may extract a further subsequence $(K_{i_k})$ along which, for each $j$, $\phi_{j,i_k}\to v_j$ weakly in $H^1$ and strongly in $L^2$ as $k\to\infty$.

If $A\subseteq RB^n\setminus K$ is compact, then $A\cap K_i=\varnothing$ for large $i$, so $\phi_{j,i}=0$ on $A$. Hence $v_j=0$ almost everywhere outside $K$. Since $\Omega$ is Lipschitz, $v_j|_\Omega\in H_0^1(\Omega)$ \cite[Corollary~1.4.4.5]{Grisvard}. For $\varphi\in C_c^\infty(\Omega)$, we have $\supp\varphi\subseteq\intt K_i$ for large $i$. Testing the eigenvalue equation for $\phi_{j,i_k}$ with $\varphi$ and letting $k\to\infty$, we obtain
\begin{equation}\label{eq:limit-weak-eigenvalue}
 \int_\Omega\nabla v_j\cdot\nabla\varphi\,dx=\lambda_j(K)\int_\Omega v_j\varphi\,dx,\qquad \varphi\in C_c^\infty(\Omega).
\end{equation}
Moreover, strong $L^2$ convergence gives $\|v_j\|_2=1$. Since $C_c^\infty(\Omega)$ is dense in $H_0^1(\Omega)$, the weak equation \eqref{eq:limit-weak-eigenvalue} also holds for $\varphi=v_j$. Hence
\[
 \|\nabla\phi_{j,i_k}\|_2^2=\lambda_j(K_{i_k})\longrightarrow\lambda_j(K)=\|\nabla v_j\|_2^2,
\]
so $\phi_{j,i_k}\to v_j$ strongly in $H^1$.
\vskip 3pt
\textbf{Step 2.} For $D\in\Ko$ and Borel $A\subseteq\partial D$, define
\[
 \gamma_D(A)=\sigma\int_A h_D(\nu_D(x))\sum_{j=1}^\infty\lambda_j(D)^{-\sigma-1}|\partial_\nu\phi_j(x)|^2\,d\mathcal{H}^{n-1}(x),
\]
where $\{\phi_j\}_{j\geqslant1}$ is an orthonormal Dirichlet eigenbasis on $\intt D$. We prove that $\gamma_{K_i}\to\gamma_K$ weakly on $RB^n$.

Let $u$ be an $L^2$-normalized Dirichlet eigenfunction on $\intt D$. For $g\in C^1(\R^n)$, write
\[
 I_D(u,g)=\int_{\partial D}g(x)(x\cdot\nu_D)|\partial_\nu u|^2\,d\mathcal{H}^{n-1}(x).
\]
Taking $X(x)=g(x)x$ in the Rellich identity \eqref{sp:eq:zetarellich} and using the strong $H^1$ convergence from Step~1, we obtain $I_{K_{i_k}}(\phi_{j,i_k},g)\to I_K(v_j,g)$.

The definition of $\gamma_D$ gives
\[
 \int_{\partial D}g\,d\gamma_D=\sigma\sum_{j=1}^\infty\lambda_j(D)^{-\sigma-1}I_D(\phi_j,g).
\]
By the Rellich identity \eqref{sp:eq:zetarellich}, $I_D(\phi_j,1)=2\lambda_j(D)$. Since $D\subseteq RB^n$ for $D=K_i$ or $K$, it follows that
\[
 \sigma\lambda_j(D)^{-\sigma-1}|I_D(\phi_j,g)|
 \leqslant 2\sigma\|g\|_{C(RB^n)}\lambda_j(D)^{-\sigma}
 \leqslant 2\sigma\|g\|_{C(RB^n)}\lambda_j(RB^n)^{-\sigma}.
\]
Since $\zeta_\sigma(RB^n)<\infty$, we may let $k\to\infty$ term by term in the sum over $j$. The $v_j$ remain orthonormal by strong $L^2$ convergence. Thus $\{v_j:\lambda_j(K)=\lambda\}$ is a basis of the $\lambda$-eigenspace. Since the sum is invariant under an orthogonal change of basis in this eigenspace, it follows that
\[
 \sum_{\lambda_j(K)=\lambda} I_K(v_j,g)=\sum_{\lambda_j(K)=\lambda} I_K(\phi_j,g).
\]
Thus $\int_{\partial K_{i_k}}g\,d\gamma_{K_{i_k}}\to\int_{\partial K}g\,d\gamma_K$.

Since the initial subsequence was arbitrary, the same limit holds for the full sequence. In particular, taking $g=1$, we have $\gamma_{K_i}(\partial K_i)\to\gamma_K(\partial K)$ as $i\to\infty$. Given $g\in C(RB^n)$ and $\varepsilon>0$, choose $h\in C^1(\R^n)$ with $\|g-h\|_{C(RB^n)}<\varepsilon$. Since the integrals of $h$ converge, the triangle inequality yields
\[
 \limsup_{i\to\infty}\Bigl|\int_{\partial K_i}g\,d\gamma_{K_i}-\int_{\partial K}g\,d\gamma_K\Bigr|\leqslant\varepsilon\lim_{i\to\infty}\bigl(\gamma_{K_i}(\partial K_i)+\gamma_K(\partial K)\bigr)=2\varepsilon\gamma_K(\partial K).
\]
Letting $\varepsilon\to0$, we conclude that $\gamma_{K_i}\to\gamma_K$ weakly on $RB^n$.

\vskip 3pt
\textbf{Step 3.} Prove that $\zeta_\sigma(K_i,\cdot)\to\zeta_\sigma(K,\cdot)$ weakly.

Let $\widehat{\gamma}_{K_i}=(\operatorname{id},\nu_{K_i})_\#\gamma_{K_i}$. By Step~2, $\sup_i\widehat{\gamma}_{K_i}(RB^n\times\Sph)=\sup_i\gamma_{K_i}(\partial K_i)<\infty$. Since $RB^n\times\Sph$ is compact, every subsequence has a weakly convergent further subsequence. Since $\gamma_{K_i}\to\gamma_K$ weakly, any such limit $\widehat{\gamma}$ satisfies
\[
 \int_{RB^n\times\Sph}f(x)\,d\widehat{\gamma}(x,u)=\int_{\partial K}f(x)\,d\gamma_K(x)\qquad\text{for }f\in C(RB^n).
\]
Since $x\cdot u=h_{K_i}(u)$ for $\widehat{\gamma}_{K_i}$-almost every $(x,u)$, it follows that
\[
 \int_{RB^n\times\Sph}|x\cdot u-h_K(u)|\,d\widehat{\gamma}_{K_i}(x,u)
 \leqslant\|h_{K_i}-h_K\|_\infty\,\gamma_{K_i}(\partial K_i)\longrightarrow0.
\]
Since the integrand is continuous, $x\cdot u=h_K(u)$ for $\widehat{\gamma}$-almost every $(x,u)$. The measure $\gamma_K$ is absolutely continuous with respect to boundary area, so $\gamma_K(\partial K\setminus\partial'K)=0$. For $x\in\partial'K$, $x\cdot u=h_K(u)$ implies $u=\nu_K(x)$. Hence, for $g\in C(\Sph)$, \eqref{sp:eq:zetameasure} gives
\[
 \int_{RB^n\times\Sph}g(u)\,d\widehat{\gamma}(x,u)=\int_{\partial K}g(\nu_K(x))\,d\gamma_K(x)=\int_{\Sph}g(u)\,d\zeta_\sigma(K,u).
\]
Since $\widehat{\gamma}$ was an arbitrary subsequential weak limit, $\zeta_\sigma(K_i,\cdot)\to\zeta_\sigma(K,\cdot)$ weakly on $\Sph$ as $i\to\infty$.
\end{proof}

\subsection{Variation formulas}\label{app:variations}

In this subsection, we verify \hyperref[cond:variation]{(V)} for $\mathcal F\in\{\lambda_1^{-\frac{1}{2}},\zeta_\sigma^{\frac{1}{2\sigma}},\tau^{\frac{1}{n+2}},I_q^{\frac{1}{n+q-1}}\}$. We begin with the variation formula for $\zeta_\sigma$.

\begin{lemma}\label{sp:prop:zetawulff}
For $K\in\Ko$, $f\in C(\Sph)$, and $\sigma>\frac{n}{2}$,
\begin{equation}\label{sp:eq:zetalinearvar}
 \frac{d}{dt}\zeta_\sigma([h_K+tf])\Bigr|_{t=0}=\int_{\Sph}\frac{f(u)}{h_K(u)}\,d\zeta_\sigma(K,u).
\end{equation}
\end{lemma}

\begin{proof}
We divide the proof into three steps.
\vskip 3pt
\textbf{Step 1.} Prove the formula for smooth bodies.

Suppose that $\partial K$ is smooth with positive Gauss curvature and $f\in C^\infty(\Sph)$. For small $|t|$, $h_K+tf$ is the support function of a body $K_t$ of the same class. The point $x_t(u)\in\partial K_t$ with outer unit normal $u$ is
\[
 x_t(u)=x_0(u)+t\bigl(\nabla_{\Sph}f(u)+f(u)u\bigr).
\]
Let $V$ be a smooth extension of the boundary field $V(x_0(u))=\nabla_{\Sph}f(u)+f(u)u$ to a neighborhood of $K$. Then $K_t=(\mathrm{id}+tV)(K)$.

For an eigenvalue $\lambda$ of $K$, taking the trace on its eigenspace in the variation formula of El Soufi--Ilias \cite[Lemmas~3.1 and~3.2]{EI} gives
\[
 \frac{d}{dt}\sum_{j:\,\lambda_j(K)=\lambda}\lambda_j(K_t)\biggr|_{t=0}
 =-\int_{\partial K}f(\nu_K)\sum_{j:\,\lambda_j(K)=\lambda}|\partial_\nu\phi_j|^2\,d\mathcal{H}^{n-1}.
\]
Choose $N$ with $\lambda_N(K)<\lambda_{N+1}(K)$. Multiplying the preceding identity by $-\sigma\lambda^{-\sigma-1}$ and summing over the distinct eigenvalues $\lambda\leqslant\lambda_N(K)$, we obtain
\begin{equation}\label{sp:eq:zeta-partial-variation}
 \frac{d}{dt}\sum_{j\leqslant N}\lambda_j(K_t)^{-\sigma}\biggr|_{t=0}=\sigma\int_{\partial K}f(\nu_K)\sum_{j\leqslant N}\lambda_j(K)^{-\sigma-1}|\partial_\nu\phi_j|^2\,d\mathcal{H}^{n-1}.
\end{equation}

Write $c=\bigl\|\frac{f}{h_K}\bigr\|_\infty$. Then $(1-c|t|)K\subseteq K_t\subseteq(1+c|t|)K$ for small $|t|$. For small $t\ne0$, it follows from \eqref{eq:eigenvalue-monotonicity-scaling} that
\[
 \sum_{j>N}\left|\frac{\lambda_j(K_t)^{-\sigma}-\lambda_j(K)^{-\sigma}}{t}\right|
 \leqslant C\sum_{j>N}\lambda_j(K)^{-\sigma}\longrightarrow0
 \qquad\text{as }N\to\infty,
\]
where $C$ is independent of $t$ and $N$, and the limit follows from $\zeta_\sigma(K)<\infty$. The series in \eqref{sp:eq:zetameasure} converges in $L^1(\partial K)$. Thus, letting $N\to\infty$ in \eqref{sp:eq:zeta-partial-variation} through indices with $\lambda_N(K)<\lambda_{N+1}(K)$ gives \eqref{sp:eq:zetalinearvar}.
\vskip 3pt
\textbf{Step 2.} Prove a mean value formula.

For $K\in\Ko$, write $d\mu_\zeta(K,u)=\frac{1}{h_K(u)}\,d\zeta_\sigma(K,u)$. Since $\frac{1}{h_{K_i}}\to\frac{1}{h_K}$ uniformly on $\Sph$ when $K_i\to K$ in $\Ko$, Lemma~\ref{sp:prop:zetaweak} implies the weak continuity of $\mu_\zeta$.

For $K,L\in\Ko$, let $M_\theta=(1-\theta)K+\theta L$. Then $h_{M_\theta}=h_K+\theta(h_L-h_K)$. Assume $K$ and $L$ are smooth with positive Gauss curvature. By Step~1 with $f=h_L-h_K$ and the mean value theorem applied to $s\mapsto\zeta_\sigma(M_s)$, there exists $\theta\in(0,1)$ such that
\begin{equation}\label{sp:eq:zeta-segment-mean-value}
 \zeta_\sigma(L)-\zeta_\sigma(K)=\int_{\Sph}(h_L(u)-h_K(u))\,d\mu_\zeta(M_\theta,u).
\end{equation}
For general $K,L$, choose smooth approximations $K_i\to K$ and $L_i\to L$ with positive Gauss curvature \cite[Section~3.4]{Schneider}. Passing to a subsequence, we may assume that the corresponding parameters $\theta_i$ tend to some $\theta\in[0,1]$. Then $(1-\theta_i)K_i+\theta_i L_i\to M_\theta$ and $h_{L_i}-h_{K_i}\to h_L-h_K$ uniformly on $\Sph$. By the continuity of $\zeta_\sigma$ and the weak continuity of $\mu_\zeta$, \eqref{sp:eq:zeta-segment-mean-value} holds for this $\theta$.

\vskip 3pt
\textbf{Step 3.} Extend the variation formula to general $K$ and continuous $f$.

Let $c=\bigl\|\frac{f}{h_K}\bigr\|_\infty$. For $t\ne0$ with $c|t|<1$, let $K_t=[h_K+tf]$ and $f_t=\frac{h_{K_t}-h_K}{t}$. Then $(1-c|t|)K\subseteq K_t\subseteq(1+c|t|)K$, and hence $\|f_t\|_\infty\leqslant c\|h_K\|_\infty$. By \eqref{sp:eq:zeta-segment-mean-value} with $L=K_t$, there exists $\theta_t\in[0,1]$ such that
\[
 \frac{\zeta_\sigma(K_t)-\zeta_\sigma(K)}{t}=\int_{\Sph}f_t(u)\,d\mu_\zeta(M_t,u),
\]
where $M_t=(1-\theta_t)K+\theta_t K_t$. Since $\|h_{M_t}-h_K\|_\infty=\theta_t|t|\,\|f_t\|_\infty\to0$, we have $M_t\to K$ as $t\to0$.

As $t\to0$, $f_t\to f$ almost everywhere with respect to $S_K$ (see, e.g., \cite[Lemma~3.3]{Jerison} and \cite[Theorem~3.5]{SSZ}). Since $\mu_\zeta(K,\cdot)$ is absolutely continuous with respect to $S_K$, the convergence also holds almost everywhere with respect to $\mu_\zeta(K,\cdot)$.

Fix $\varepsilon>0$ with $c\varepsilon<1$. With $K_0=K$, the definition of $K_t$ gives $(1-\theta)K_{t_0}+\theta K_{t_1}\subseteq K_{(1-\theta)t_0+\theta t_1}$ for $t_0,t_1\in[-\varepsilon,\varepsilon]$ and $\theta\in[0,1]$. Thus $t\mapsto h_{K_t}(u)$ is concave and $f_\varepsilon\leqslant f_t\leqslant f_{-\varepsilon}$ for $0<|t|\leqslant\varepsilon$. Integrating against $\mu_\zeta(M_t,\cdot)$, we obtain
\[
 \int_{\Sph}f_\varepsilon(u)\,d\mu_\zeta(M_t,u)
 \leqslant\frac{\zeta_\sigma(K_t)-\zeta_\sigma(K)}{t}
 \leqslant\int_{\Sph}f_{-\varepsilon}(u)\,d\mu_\zeta(M_t,u).
\]
Since $f_{\pm\varepsilon}$ are continuous and $\mu_\zeta(M_t,\cdot)$ converges weakly to $\mu_\zeta(K,\cdot)$ as $t\to0$, we have
\[
 \lim_{t\to0}\int_{\Sph}f_{\pm\varepsilon}(u)\,d\mu_\zeta(M_t,u)=\int_{\Sph}f_{\pm\varepsilon}(u)\,d\mu_\zeta(K,u).
\]
As $\varepsilon\to0^+$, the dominated convergence theorem gives the common limit $\int_{\Sph}f(u)\,d\mu_\zeta(K,u)$. Therefore, \eqref{sp:eq:zetalinearvar} holds.
\end{proof}

We now turn to the measures associated with the other three functionals.

Following Crasta and Fragal\`a \cite[Definition~1]{CrF}, define the cone eigenvalue measure and the cone torsion measure of $K\in\Ko$ by
\begin{align*}
 \lambda_1(K,\omega)&=\int_{\nu_K^{-1}(\omega)}h_K(\nu_K(x))|\partial_\nu\phi_1(x)|^2\,d\mathcal{H}^{n-1}(x),\\
 \tau(K,\omega)&=\int_{\nu_K^{-1}(\omega)}h_K(\nu_K(x))|\nabla u_K(x)|^2\,d\mathcal{H}^{n-1}(x),
\end{align*}
for Borel $\omega\subseteq\Sph$. Here, the boundary values of $\nabla u_K$ are taken as nontangential limits.

Following Lutwak, Xi, Yang, and Zhang \cite[Definition~4.1 and (2.11)]{LXYZ}, define the $q$th cone-chord measure of $K$ for $1<q<n+1$ by
\[
 G_q(K,\omega)=\frac{2q(q-1)}{n\omega_n(n+q-1)}\int_{\nu_K^{-1}(\omega)}h_K(\nu_K(x))\int_K |x-y|^{q-1-n}\,dy\,d\mathcal{H}^{n-1}(x).
\]

To verify \hyperref[cond:variation]{(V)}, let $K\in\Ko$, $g\in C(\Sph)$. Write $K_t=[h_Ke^{tg}]$ for $t\in\mathbb{R}$. For $\mathcal F\in\{\lambda_1^{-\frac{1}{2}},\zeta_\sigma^{\frac{1}{2\sigma}},\tau^{\frac{1}{n+2}},I_q^{\frac{1}{n+q-1}}\}$ and $|t|$ small, Lemma~\ref{lem:wulff-comparison} gives
\[
 \bigl|\mathcal F(K_t)-\mathcal F([h_K+t h_Kg])\bigr|\leqslant\mathcal F([h_K+t h_Kg])\left\|\frac{e^{tg}}{1+tg}-1\right\|_\infty=O(t^2).
\]
Hence $\frac{d}{dt}\mathcal F(K_t)\bigr|_{t=0}=\frac{d}{dt}\mathcal F([h_K+t h_Kg])\bigr|_{t=0}$ whenever the derivative on the right exists. Therefore, taking $f=h_Kg$ in the formula of Jerison \cite[Theorem~7.5]{Jerison}, Lemma~\ref{sp:prop:zetawulff}, and the formula of Colesanti--Fimiani \cite[Corollary~4.2]{CF} gives the first three variational formulas below by the chain rule. The fourth follows from the formula of Lutwak--Xi--Yang--Zhang \cite[equation~(5.11)]{LXYZ}.
\begin{equation}\label{eq:normalized-log-variations}
\begin{aligned}
 \frac{d}{dt}\log\bigl(\lambda_1(K_t)^{-\frac{1}{2}}\bigr)\Bigr|_{t=0}
 &=\frac{1}{2\lambda_1(K)}\int_{\Sph}g\,d\lambda_1(K,u),\\
 \frac{d}{dt}\log\bigl(\zeta_\sigma(K_t)^{\frac{1}{2\sigma}}\bigr)\Bigr|_{t=0}
 &=\frac{1}{2\sigma\zeta_\sigma(K)}\int_{\Sph}g\,d\zeta_\sigma(K,u),\\
 \frac{d}{dt}\log\bigl(\tau(K_t)^{\frac{1}{n+2}}\bigr)\Bigr|_{t=0}
 &=\frac{1}{(n+2)\tau(K)}\int_{\Sph}g\,d\tau(K,u),\\
 \frac{d}{dt}\log\bigl(I_q(K_t)^{\frac{1}{n+q-1}}\bigr)\Bigr|_{t=0}
 &=\frac{1}{I_q(K)}\int_{\Sph}g\,dG_q(K,u).
\end{aligned}
\end{equation}
If $K\in\Ke$ and $g$ is odd, then $K_{-t}=-K_t$, so $\mathcal F(K_{-t})=\mathcal F(K_t)$. Thus the derivatives in \eqref{eq:normalized-log-variations} vanish for odd $g$, and hence the four measures are even. This verifies \hyperref[cond:variation]{(V)} for the four normalized functionals.

\subsection{Estimates for geometric functionals}\label{app:growth}

In this subsection, we verify \hyperref[cond:john]{(J)} and \hyperref[cond:cylinder]{(C)} for $\mathcal F\in\{\lambda_1^{-\frac{1}{2}},\zeta_\sigma^{\frac{1}{2\sigma}},\tau^{\frac{1}{n+2}},I_q^{\frac{1}{n+q-1}}\}$.

For $K\in\Ke$, let $0<a_1\leqslant\cdots\leqslant a_n$ be the semiaxis lengths of its John ellipsoid $\mathrm{J}K$. In principal coordinates, write $Q=\prod_{i=1}^n[-a_i,a_i]$. By John's theorem (see, e.g., \cite[Theorem~10.12.2]{Schneider}), we have
\[
 \frac{1}{\sqrt n}Q\subseteq\mathrm{J}K\subseteq K\subseteq\sqrt n\,\mathrm{J}K\subseteq\sqrt n\,Q.
\]
Since $\mathcal F$ is monotone and homogeneous of degree one, these inclusions give $\mathcal F(K)\asymp\mathcal F(Q)$. Here $f\asymp g$ means that $C^{-1}g\leqslant f\leqslant Cg$ for some $C\geqslant1$ independent of the half-side lengths $a_1,\ldots,a_n$ of $Q$.

\begin{lemma}\label{sp:prop:john-bounds}
For $\sigma>\frac{n}{2}$,
\[
 \lambda_1(Q)\asymp a_1^{-2},\qquad \zeta_\sigma(Q)\asymp |Q|a_1^{2\sigma-n},\qquad \tau(Q)\asymp |Q|a_1^2.
\]
\end{lemma}
\begin{proof}
By separation of variables, $\lambda_1(Q)$ is the sum of the first Dirichlet eigenvalues of the intervals $(-a_i,a_i)$. Thus $\lambda_1(Q)=\sum_{i=1}^n\frac{\pi^2}{4a_i^2}\asymp a_1^{-2}$, since $a_1=\min_i a_i$.

For $t>0$, \eqref{eq:heat-trace-expansion} and \eqref{eq:heat-trace-upper-bound} give
\[
 Z_Q(t)\leqslant e^{-\frac{t}{2}\lambda_1(Q)}Z_Q\Bigl(\frac{t}{2}\Bigr)\leqslant |Q|(2\pi t)^{-\frac{n}{2}}e^{-\frac{t}{2}\lambda_1(Q)}.
\]
It follows from \eqref{eq:zeta-heat-representation} that
\[
 \zeta_\sigma(Q)\leqslant \frac{|Q|}{(2\pi)^{\frac{n}{2}}\Gamma(\sigma)}\int_0^\infty t^{\sigma-\frac{n}{2}-1}e^{-\frac{t}{2}\lambda_1(Q)}\,dt=C_{n,\sigma}|Q|\lambda_1(Q)^{\frac{n}{2}-\sigma}\asymp |Q|a_1^{2\sigma-n}.
\]

For torsional rigidity, by \eqref{eq:torsion-spectral-expansion} and $\lambda_j(Q)\geqslant\lambda_1(Q)$, we have
\[
 \tau(Q)=\sum_{j=1}^\infty\frac{\bigl(\int_Q\phi_j\,dx\bigr)^2}{\lambda_j(Q)}
 \leqslant\frac{1}{\lambda_1(Q)}\sum_{j=1}^\infty\Bigl(\int_Q\phi_j\,dx\Bigr)^2
 =\frac{|Q|}{\lambda_1(Q)}\asymp |Q|a_1^2.
\]
The last equality follows from Parseval's identity
\[
 \sum_{j=1}^\infty\Bigl(\int_Q\phi_j\,dx\Bigr)^2
 =\sum_{j=1}^\infty\bigl|\langle 1,\phi_j\rangle_{L^2(Q)}\bigr|^2
 =\|1\|_{L^2(Q)}^2=|Q|.
\]

For the lower bounds, let $U=[-a_1,a_1]^n$. Then $Q$ contains $N=\prod_{i=1}^n\bigl\lfloor\frac{a_i}{a_1}\bigr\rfloor\asymp |Q|a_1^{-n}$ translates of $U$ with disjoint interiors. For each such translate $D$ and $t>0$, domain monotonicity of the Dirichlet heat kernel (see Section~\ref{sec:heat-kernels}) implies
\begin{align*}
 \int_D k_t^Q(x,x)\,dx &\geqslant \int_D k_t^D(x,x)\,dx=Z_D(t)=Z_U(t),\\
 \int_{D^2}k_t^Q(x,y)\,dx\,dy &\geqslant \int_{D^2}k_t^D(x,y)\,dx\,dy=Q_D(t)=Q_U(t).
\end{align*}
Summing over the $N$ cubes and applying \eqref{eq:zeta-heat-representation} and \eqref{eq:torsion-heat-representation}, we obtain
\[
 \zeta_\sigma(Q)\geqslant N\zeta_\sigma(U)\asymp |Q|a_1^{2\sigma-n},\qquad
 \tau(Q)\geqslant N\tau(U)\asymp |Q|a_1^2,
\]
where $\zeta_\sigma(U)\asymp a_1^{2\sigma}$ and $\tau(U)\asymp a_1^{n+2}$ by homogeneity. This completes the proof.
\end{proof}

The estimates for $I_q(Q)$ in the next lemma follow from the dual volume estimates of B\"or\"oczky--Kov\'acs--Mui--Zhang \cite[Lemma~4.1]{BKMZ}. For noninteger $q$, the ellipsoid estimate of Lutwak--Xi--Yang--Zhang \cite[Lemma~7.3]{LXYZ} already gives the upper bound. We give a direct proof here.
\begin{lemma}\label{prop:riesz-box-bounds}
Let $1<q<n+1$ and $\ell=\lfloor q\rfloor$. If $\ell<q<\ell+1$, then
\begin{equation}\label{ri:eq:box-noninteger}
 I_q(Q)\asymp |Q|\Bigl(\prod_{i<\ell}a_i\Bigr)a_\ell^{q-\ell}.
\end{equation}
If $q=\ell$, then
\begin{equation}\label{ri:eq:box-integer}
 I_q(Q)\asymp |Q|\Bigl(\prod_{i<\ell}a_i\Bigr)\Bigl(1+\log\frac{a_\ell}{a_{\ell-1}}\Bigr).
\end{equation}
\end{lemma}
\begin{proof}
By \eqref{ri:eq:chord}, it suffices to estimate $E_{n+1-q}(Q)$. Using the definition in Section~\ref{sec:geometric-functionals} and Fubini's theorem, we have
\begin{align*}
 E_{n+1-q}(Q)
 &=\int_Q\int_Q |x-y|^{q-n-1}\,dy\,dx\\
 &=(n+1-q)\int_Q\int_Q\int_{|x-y|}^\infty r^{q-n-2}\,dr\,dy\,dx\\
 &=(n+1-q)\int_0^\infty r^{q-n-2}\int_Q\int_{Q\cap B(x,r)}dy\,dx\,dr\\
 &=(n+1-q)\int_0^\infty r^{q-n-2}\int_Q |Q\cap B(x,r)|\,dx\,dr,
\end{align*}
where $B(x,r)=x+rB^n$.

To estimate $|Q\cap B(x,r)|$, note that
\[
 x+\frac{r}{\sqrt n}[-1,1]^n\subseteq B(x,r)\subseteq x+r[-1,1]^n.
\]
For $x\in Q$ and $s>0$, we have $\min\{a_i,s\}\leqslant\bigl|[-a_i,a_i]\cap[x_i-s,x_i+s]\bigr|\leqslant 2\min\{a_i,s\}$. So
\[
 \bigl|Q\cap(x+[-s,s]^n)\bigr|=\prod_{i=1}^n\bigl|[-a_i,a_i]\cap[x_i-s,x_i+s]\bigr|\asymp\prod_{i=1}^n\min\{a_i,s\}.
\]
Taking $s=\frac{r}{\sqrt n}$ and $s=r$ yields $|Q\cap B(x,r)|\asymp\prod_{i=1}^n\min\{a_i,r\}$, uniformly for $x\in Q$ and $r>0$. Hence
\begin{equation}\label{ri:eq:rectangle-integral}
 I_q(Q)\asymp |Q|\int_0^\infty H(r)\,\frac{dr}{r},\qquad H(r)=r^{q-n-1}\prod_{i=1}^n\min\{a_i,r\}.
\end{equation}

Write $a_0=0$, $a_{n+1}=\infty$, and define $H(0)=H(\infty)=0$. For $a_{j-1}<r<a_j$, $1\leqslant j\leqslant n+1$, we have $H(r)=\bigl(\prod_{i<j}a_i\bigr)r^{q-j}$. Thus $H$ is nondecreasing on $[0,a_\ell]$ and nonincreasing on $[a_\ell,\infty)$. Its maximum is $M=H(a_\ell)=\bigl(\prod_{i<\ell}a_i\bigr)a_\ell^{q-\ell}$. For $j\ne q$,
\[
 \int_{a_{j-1}}^{a_j}H(r)\,\frac{dr}{r}=\frac{|H(a_j)-H(a_{j-1})|}{|q-j|}\asymp |H(a_j)-H(a_{j-1})|.
\]
Summing over $j\ne q$, we obtain
\[
 \sum_{j\ne q}\int_{a_{j-1}}^{a_j}H(r)\,\frac{dr}{r}\asymp\sum_{j\ne q}|H(a_j)-H(a_{j-1})|=2M.
\]

If $\ell<q<\ell+1$, then
\[
 \int_0^\infty H(r)\,\frac{dr}{r}=\sum_{j\ne q}\int_{a_{j-1}}^{a_j}H(r)\,\frac{dr}{r}\asymp M.
\]
If $q=\ell$, then $H=M$ on $[a_{\ell-1},a_\ell]$. Hence
\[
 \int_0^\infty H(r)\,\frac{dr}{r}\asymp M+M\int_{a_{\ell-1}}^{a_\ell}\frac{dr}{r}=M\Bigl(1+\log\frac{a_\ell}{a_{\ell-1}}\Bigr).
\]
Substituting these estimates into \eqref{ri:eq:rectangle-integral} yields \eqref{ri:eq:box-noninteger} and \eqref{ri:eq:box-integer}, respectively.
\end{proof}

We now verify \hyperref[cond:john]{(J)} and \hyperref[cond:cylinder]{(C)}. For the four normalized functionals, Table~\ref{tab:inverse-inputs} lists the exponents $b_i$ in \hyperref[cond:john]{(J)} obtained from Lemmas~\ref{sp:prop:john-bounds} and~\ref{prop:riesz-box-bounds}.

\begin{table}[ht]
\centering
\renewcommand{\arraystretch}{1.25}
\setlength{\tabcolsep}{9pt}
\begin{tabular}{lc@{\hspace{24pt}}c}
\toprule
$\mathcal F$ & $(b_1,\ldots,b_n)$ & $c_k=\sum_{i=1}^k b_i$ \\
\midrule
$\lambda_1^{-\frac{1}{2}}$ & $(1,0,\ldots,0)$ & $1$ \\\addlinespace[3pt]
$\zeta_\sigma^{\frac{1}{2\sigma}}$ & $\dfrac{(2\sigma-n+1,1,\ldots,1)}{2\sigma}$ & $1-\dfrac{n-k}{2\sigma}$ \\\addlinespace[3pt]
$\tau^{\frac{1}{n+2}}$ & $\dfrac{(3,1,\ldots,1)}{n+2}$ & $\dfrac{k+2}{n+2}$ \\\addlinespace[3pt]
$I_q^{\frac{1}{n+q-1}}$ & $\dfrac{(1+\min\{i,q-1\}-\min\{i-1,q-1\})_{i=1}^n}{n+q-1}$ & $\dfrac{k+\min\{k,q-1\}}{n+q-1}$ \\
\bottomrule
\end{tabular}
\caption{}
\label{tab:inverse-inputs}
\end{table}

For $\mathcal F\in\{\lambda_1^{-\frac{1}{2}},\zeta_\sigma^{\frac{1}{2\sigma}},\tau^{\frac{1}{n+2}}\}$, Lemma~\ref{sp:prop:john-bounds} gives $\mathcal F(Q)\asymp\prod_{i=1}^n a_i^{b_i}$. By Lemma~\ref{prop:riesz-box-bounds}, the same estimate holds for $\mathcal F=I_q^{\frac{1}{n+q-1}}$ when $q\notin\mathbb Z$. If $q=\ell\in\mathbb Z$, then
\[
 \mathcal F(Q)\asymp\Bigl(\prod_{i=1}^n a_i^{b_i}\Bigr)\Bigl(1+\log\frac{a_\ell}{a_{\ell-1}}\Bigr)^{\frac{1}{n+q-1}}.
\]
Observe that for any $\varepsilon>0$, $\bigl(1+\log\frac{a_\ell}{a_{\ell-1}}\bigr)^{\frac{1}{n+q-1}}\leqslant A_\varepsilon\bigl(\frac{a_n}{a_1}\bigr)^\varepsilon$, where $A_\varepsilon>0$ is independent of $a_1,\ldots,a_n$. Since $\mathcal F(K)\asymp\mathcal F(Q)$, all four functionals satisfy
\[
 \mathcal F(K)\leqslant C_\varepsilon\Bigl(\frac{a_n}{a_1}\Bigr)^\varepsilon\prod_{i=1}^n a_i^{b_i},
\]
as required in \hyperref[cond:john]{(J)}.

To verify \hyperref[cond:cylinder]{(C)}, fix a $k$-dimensional subspace $\xi$ of $\R^n$ and choose orthonormal coordinates so that $\xi=\R^k\times\{0\}$. For $r\geqslant 1$, let $Q_r=[-1,1]^k\times[-r,r]^{n-k}$. Since $\frac{1}{\sqrt n}Q_r\subseteq B_\xi\times rB_{\xi^\perp}\subseteq Q_r$, we have $\mathcal F(B_\xi\times rB_{\xi^\perp})\asymp\mathcal F(Q_r)$.

Since $\mathrm{J}Q_r$ has semiaxis lengths $a_1=\cdots=a_k=1$ and $a_{k+1}=\cdots=a_n=r$, we have
\[
 \prod_{i=1}^n a_i^{b_i}=r^{\sum_{i=k+1}^n b_i}=r^{1-c_k}.
\]
Applying \hyperref[cond:john]{(J)} to $Q_r$ yields $\mathcal F(Q_r)\leqslant C_\varepsilon r^{1-c_k+\varepsilon}$ for $\varepsilon>0$. Lemmas~\ref{sp:prop:john-bounds} and~\ref{prop:riesz-box-bounds} give $\mathcal F(Q_r)\geqslant cr^{1-c_k}$. Hence
\[
 1-c_k\leqslant\liminf_{r\to\infty}\frac{\log\mathcal F(Q_r)}{\log r}\leqslant\limsup_{r\to\infty}\frac{\log\mathcal F(Q_r)}{\log r}\leqslant 1-c_k+\varepsilon.
\]
Letting $\varepsilon\to0$ and using $\mathcal F(B_\xi\times rB_{\xi^\perp})\asymp\mathcal F(Q_r)$, we conclude that
\[
 \lim_{r\to\infty}\frac{\log\mathcal F(B_\xi\times rB_{\xi^\perp})}{\log r}=\lim_{r\to\infty}\frac{\log\mathcal F(Q_r)}{\log r}=1-c_k,
\]
as required in \hyperref[cond:cylinder]{(C)}.

\subsection{Proof of Theorem~\ref{thm:geometric-inverse}}\label{sec:geometric-inverse-proof}

\begin{proof}
In (i)--(iv), take $(\mathcal F,\widetilde{\eta}_K)$ to be
\[
 \bigl(\lambda_1^{-\frac{1}{2}},\lambda_1(K,\cdot)\bigr),\quad
 \bigl(\zeta_\sigma^{\frac{1}{2\sigma}},\zeta_\sigma(K,\cdot)\bigr),\quad
 \bigl(\tau^{\frac{1}{n+2}},\tau(K,\cdot)\bigr),\quad
 \bigl(I_q^{\frac{1}{n+q-1}},G_q(K,\cdot)\bigr),
\]
respectively. Then $\mathcal F$ satisfies \hyperref[cond:logbm]{(L)} by Theorem~\ref{thm:geometric-logbm} and \hyperref[cond:john]{(J)} and \hyperref[cond:cylinder]{(C)} by Section~\ref{app:growth}. By \eqref{eq:normalized-log-variations}, $\mathcal F$ satisfies \hyperref[cond:variation]{(V)} with $\eta_K$ a positive multiple of $\widetilde{\eta}_K$. Since $\eta_K$ is a probability measure, we have $\eta_K=\frac{\widetilde{\eta}_K}{\widetilde{\eta}_K(\Sph)}$.

Assume that $\widetilde{\eta}_K=\mu$ for some $K\in\Ke$. Fix a nonzero proper subspace $\xi$ of $\R^n$ and write $k=\dim\xi$. By Theorem~\ref{thm:abstract-inverse-principle}~(i), $\eta_K\in\mathcal P_c$, and hence
\[
 \frac{\mu(\Sph\cap\xi)}{\mu(\Sph)}=\frac{\widetilde{\eta}_K(\Sph\cap\xi)}{\widetilde{\eta}_K(\Sph)}=\eta_K(\Sph\cap\xi)<c_k.
\]
Substituting the values of $c_k$ from Table~\ref{tab:inverse-inputs} proves necessity in (i)--(iv).

Conversely, if $\mu$ satisfies the corresponding subspace concentration condition, then $\frac{\mu}{\mu(\Sph)}\in\mathcal P_c$. By Theorem~\ref{thm:abstract-inverse-principle}~(ii), there exists $K\in\Ke$ such that $\eta_K=\frac{\mu}{\mu(\Sph)}$.

Taking $g=1$ in \eqref{eq:normalized-log-variations} and using the homogeneity of $\lambda_1$, $\zeta_\sigma$, $\tau$, and $I_q$, we obtain
\begin{equation}\label{eq:cone-measure-mass-scaling}
 \widetilde{\eta}_{tK}(\Sph)=t^s\widetilde{\eta}_K(\Sph),\qquad\text{for } t>0,
\end{equation}
where $s=-2,\,2\sigma,\,n+2,\,n+q-1$ in (i)--(iv), respectively. Let $r=\bigl(\frac{\mu(\Sph)}{\widetilde{\eta}_K(\Sph)}\bigr)^{\frac{1}{s}}$. Since $\eta_{rK}=\eta_K$, we have
\[
 \widetilde{\eta}_{rK}=r^s\widetilde{\eta}_K(\Sph)\,\eta_K=\mu(\Sph)\,\eta_K=\mu.
\]
This proves sufficiency.

If $L\in\Ke$ also satisfies $\widetilde{\eta}_L=\mu$, then $\eta_L=\frac{\widetilde{\eta}_L}{\widetilde{\eta}_L(\Sph)}=\frac{\mu}{\mu(\Sph)}=\eta_K$. By Theorem~\ref{thm:abstract-uniqueness}, $L=tK$ for some $t>0$. Since $\widetilde{\eta}_{tK}=\mu=\widetilde{\eta}_{rK}$, \eqref{eq:cone-measure-mass-scaling} gives $t^s=r^s$. Hence $t=r$ and $L=rK$. This proves uniqueness.
\end{proof}

\begin{remark}\label{rem:sharpness-existence}
The inequalities in Theorem~\ref{thm:geometric-inverse} are sharp for each fixed nonzero proper subspace $\xi$ of $\R^n$. Fix one of the functionals in Table~\ref{tab:inverse-inputs} and write $k=\dim\xi$. By Theorem~\ref{thm:abstract-inverse-principle}~(i), $\eta_K(\Sph\cap\xi)<c_k$ for $K\in\Ke$.

To prove sharpness, it suffices to show that for each $t\in[0,c_k)$ there exists $K_t\in\Ke$ with $\eta_{K_t}(\Sph\cap\xi)=t$. Denote by $\nu_\xi$ and $\nu$ the uniform probability measures on $\Sph\cap\xi$ and $\Sph$, respectively. Let $\mu_t=t\nu_\xi+(1-t)\nu$. For any nonzero proper subspace $\xi'$ of $\R^n$,
\[
 \mu_t(\Sph\cap \xi')=t\,\nu_\xi(\Sph\cap \xi')=
 \begin{cases}
 t,&\xi\subseteq \xi',\\
 0,&\text{otherwise}.
 \end{cases}
\]
Since $0<c_1\leqslant\cdots\leqslant c_{n-1}$, it follows that $\mu_t(\Sph\cap \xi')<c_{\dim \xi'}$. Thus $\mu_t\in\mathcal P_c$, and Theorem~\ref{thm:abstract-inverse-principle}~(ii) gives $K_t\in\Ke$ with $\eta_{K_t}=\mu_t$. Taking $\xi'=\xi$, we have $\eta_{K_t}(\Sph\cap\xi)=\mu_t(\Sph\cap\xi)=t$, as required.
\end{remark}

\vskip 25pt
\bibliographystyle{amsplain}

\end{document}